# Variance Reduction for Faster Non-Convex Optimization

Zeyuan Allen-Zhu
zeyuan@csail.mit.edu
Princeton University

Elad Hazan
ehazan@cs.princeton.edu
Princeton University

February 5, 2016*

**Abstract**

We consider the fundamental problem in non-convex optimization of efficiently reaching a stationary point. In contrast to the convex case, in the long history of this basic problem, the only known theoretical results on first-order non-convex optimization remain to be full gradient descent that converges in $O(1/\varepsilon)$ iterations for smooth objectives, and stochastic gradient descent that converges in $O(1/\varepsilon^2)$ iterations for objectives that are sum of smooth functions.

We provide the first improvement in this line of research. Our result is based on the variance reduction trick recently introduced to convex optimization, as well as a brand new analysis of variance reduction that is suitable for non-convex optimization. For objectives that are sum of smooth functions, our first-order minibatch stochastic method converges with an $O(1/\varepsilon)$ rate, and is faster than full gradient descent by $\Omega(n^{1/3})$.

We demonstrate the effectiveness of our methods on empirical risk minimizations with non-convex loss functions and training neural nets.

## 1 Introduction

Numerous machine learning problems are naturally formulated as non-convex optimization problems. Examples include inference in graphical models, unsupervised learning models such as topic models, dictionary learning, and perhaps most notably, training of deep neural networks. Indeed, non-convex optimization for machine learning is one of the fields' main research frontiers.

Since global minimization of non-convex functions is NP-hard, various alternative approaches are applied. For some models, probabilistic and other assumptions on the input can be used to give specially designed polynomial-time algorithms [5, 6, 15].

However, the multitude and diversity of machine learning applications require a robust, generic optimization method that can be applied as a tool rather than reinvented per each specific model. One approach is the design of global non-convex heuristics such as simulated annealing or bayesian optimization. Although believed to fail in the worst case due to known complexity results, such heuristics many times perform well in practice for certain problems.

Another approach, which is based on more solid theoretical foundation and is gaining in popularity, is to drop the "global optimality" requirement and attempt to reach more modest solution concepts. The most popular of these is the use of iterative optimization methods to reach a stationary point. The use of stochastic first-order methods is the primary focus of this approach, which has become the most common method for training deep neural nets.

*First circulated on this date, and first appeared on arXiv on March 17, 2016.

Formally, in this paper we consider the unconstrained minimization problem

$$\min_{x\in\mathbb{R}^d}\Big\{f(x)\stackrel{\mathrm{def}}{=}\frac{1}{n}\sum_{i=1}^{n}f_i(x)\Big\}\ , \tag{1.1}$$

where each $f_i(x)$ is differentiable, possibly non-convex, and has $L$-Lipschitz continuous gradient (a.k.a. $L$-smooth) for some parameter $L>0$.[1] Many machine learning/imaging processing problems fall into Problem (1.1), including training neural nets, ERM (empirical risk minimization) with non-convex losses, and many others.

Following the classical benchmark for non-convex optimization (see for instance [13]), we focus on algorithms that can efficiently find an approximate stationary point $x$ satisfying $\|\nabla f(x)\|^2\le\varepsilon$.

Unlike convex optimization, a point with small gradient may only be close to a saddle point or a local minimum, rather than the global minimum. Therefore, such an algorithm is usually combined with saddle-point or local-minima escaping schemes, such as genetic algorithms or simulated annealing. More recently, Ge et al. [12] also demonstrated that a simple noise-addition scheme is sufficient for stochastic gradient descent to escape from saddle points.

However, for the general problem (1.1) where smoothness is the only assumption and finding approximate stationary point is the simple goal, the only known theoretical convergence results remain to be that for *gradient descent (GD)* and *stochastic gradient descent (SGD)*.

- Given a starting point $x_0$, GD applies an update $x'\leftarrow x-\frac{1}{L}\nabla f(x)$ with a fixed step length $1/L$ per iteration. In order to produce an output $x$ that is an $\varepsilon$-approximate stationary point, GD needs $T=O\big(\frac{L(f(x_0)-f(x^*))}{\varepsilon}\big)$ iterations where $x^*$ is the global minimizer of $f(\cdot)$. This is a folklore result in optimization and included for instance in [13].
- SGD applies an update $x'\leftarrow x-\eta\nabla f_i(x)$ per iteration, where $i$ chosen uniformly at random from $[n]\stackrel{\mathrm{def}}{=}\{1,2,\dots,n\}$. If $\eta$ is properly tuned, one can obtain an $\varepsilon$-approximate stationary point in $T=O\big(\big(\frac{L}{\varepsilon}+\frac{L\sigma^2}{\varepsilon^2}\big)\cdot(f(x_0)-f(x^*))\big)$ iterations, where $\sigma$ is the variance of the stochastic gradient. This result is perhaps first formalized by Ghadimi and Lan [13].

Since computing the full gradient $\nabla f(\cdot)$ is usually $n$ times slower than that of $\nabla f_i(x)$, each iteration of SGD is usually $n$ times faster than that of GD, but the total number of iterations for SGD is very poor.

Before our work, it is an open question to design a first-order method that is faster than both GD and SGD.

## 1.1 Our Result

We prove that variance reduction techniques, based on the SVRG method [16], produce an $\varepsilon$-stationary point in only $O\big(\frac{n^{2/3}L(f(x_0)-f(x^*))}{\varepsilon}\big)$ iterations. Since each iteration of SVRG is as fast as SGD and $n$ times faster than that of GD, SVRG is guaranteed to be at least $\Omega(n^{1/3})$ times faster than GD. Among first-order methods, this is the *first time* the performance of GD is outperformed in theory for problem (1.1) without any additional assumption, and also the *first time* that stochastic-gradient based methods are shown to have a non-trivial[2] $1/\varepsilon$ convergence rate independent of the variance $\sigma^2$.

[1]Even if each $f_i(x)$ is not smooth but only Lipschitz continuous, standard smoothing techniques such as Chapter 2.3 of [14] usually turn each $f_i(x)$ into a smooth function without sacrificing too much accuracy.

[2]Note however, designing a stochastic-gradient method with a trivial $1/\varepsilon$ rate is obvious. For instance, it is straightforward to design such a method that converges in $O\big(\frac{nL(f(x_0)-f(x^*))}{\varepsilon}\big)$ iterations. However, this is never faster than GD.

Our proposed algorithm is very analogous to SVRG of [16]. Recall that SVRG has an outer loop of epochs, where at the beginning of each epoch, SVRG defines a *snapshot vector* $\widetilde{x}$ to be the average vector of the previous epoch,[3] and computes its full gradient $\nabla f(\widetilde{x})$. Each epoch of SVRG consists of $m$ inner iterations, where the choice of $m$ usually depends on the objective's strong convexity. In each inner iteration inside an epoch, SVRG picks a random $i \in [n]$, defines the *gradient estimator*

$$\widetilde{\nabla}^k \stackrel{\text{def}}{=} \frac{1}{n}\sum_{j=1}^{n} \nabla f_j(\widetilde{x}) + \nabla f_i(x^k) - \nabla f_i(\widetilde{x}) \enspace , \tag{1.2}$$

and performs an update $x' \leftarrow x - \eta\widetilde{\nabla}^k$ for some fixed step length $\eta > 0$ across all iterations and epochs.

In order to prove our theoretical result in this paper, we make the following changes to SVRG. First, we set the number of inner iterations $m$ as a constant factor times $n$. Second, we pick the snapshot point to be a non-uniform average of the last $m^{2/3}$ elements of the previous epoch. Finally, we prove that the average norm $\|\nabla f(x_k)\|^2$ of the encountered vectors $x_k$ across all iterations is small, so it suffices to output $x_k$ for a random $k$.

**Our Technique.** To prove our result, we need different techniques from all known results on variance reduction. The key idea used by previous authors is to show that the variance of the gradient estimator $\widetilde{\nabla}^k$ is upper bounded by either $O(f(x^k) - f(x^*))$ or $O(\|x^k - x^*\|^2)$, and therefore it converges to zero for convex functions. This analysis fails to apply in the non-convex setting because gradient-based methods do not converge to the global minimum.

We observe in this paper that the variance is upper bounded by $O(\|x^k - \widetilde{x}\|^2)$, the squared distance between the current point and the most recent snapshot. By dividing an epoch into $m^{1/3}$ subepochs of length $m^{2/3}$, and performing a mirror-descent analysis for each subepoch, we further show that this squared distance is related to the objective decrease $f(\widetilde{x}) - f(x^k)$. This would suffice for proving our theorem: whenever this squared distance is small the objective is decreased by a lot due to the small variance, or otherwise if this squared distance is large we still experience a large objective decrease because it is related to $f(\widetilde{x}) - f(x^k)$.

**Applications.** There are many machine learning problems that fall into category (1.1). To mention just two:

- Non-Convex Loss in ERM

  Empirical risk minimization (ERM) problems naturally fall into the category of (1.1) if the loss functions are *non-convex*. For instance, for binary classification problems, the sigmoid function —or more broadly, any natural smoothed variant of the 0-1 loss function— is not only a more natural choice than artificial ones such as hinge loss, logistic loss, squared loss, but also generalize better in terms of testing accuracy especially when there are outliers [25].

  However, since sigmoid loss is non-convex, it was previously considered hard to train an ERM problem with it. Shalev-Shwartz, Shamir and Sridharan [25] showed that this minimization problem is still solvable in the improper learning sense, with the help from kernel methods and gradient descent. However, their theoretical convergence has a poor polynomial dependence on $1/\varepsilon$ and exponential dependence on the smoothness parameter of the loss function.

  Our result in this paper applies to ERM problems with non-convex loss. Suppose we are given $n$ training examples $\{(a_1, \ell_1), \dots (a_n, \ell_n)\}$, where each $a_i \in \mathbb{R}^d$ is the feature vector of example $i$ and each $l_i \in \{-1, +1\}$ is the binary label of example $i$. By setting $\phi(t) \stackrel{\text{def}}{=} \frac{1}{1+e^t}$ to be the sigmoid loss function and setting $f_i(x) \stackrel{\text{def}}{=} \phi(l_i \langle a_i, x\rangle) + \frac{\lambda}{2}\|x\|^2$, problem (1.1) becomes $\ell_2$ ERM with sigmoid loss. We shall demonstrate in our experiment section that, by using SVRG to train ERM with sigmoid loss, its running time is as good as using SVRG to train ERM with other convex loss functions, but the testing accuracy can be significantly better.

[3] More precisely, SVRG provides two options, one defining $\widetilde{x}$ to be the average vector of the previous epoch, and the other defining $\widetilde{x}$ to be the last iterate of the previous epoch. While the authors only prove theoretical results for the "average" definition, experimental results suggest that choosing the last iterate is better.

- Neural Network

  Training neural nets can also be formalized into problem (1.1). For instance, as long as the activation function of each neural node is smooth, say the sigmoid function or a smooth version of the rectified linear unit (ReLU) function (for instance, the softplus alternative), we can define $f_i(x)$ to be the training loss with respect to the $i$-th data input. In this language, computing the stochastic gradient $\nabla f_i(x)$ for some random $i \in [n]$ corresponds to performing one forward-backward prorogation on the neural net with respect to sample $i$. We shall demonstrate in our experiment that using SVRG to train neural nets can enjoy a much faster running time comparing to SGD or SVRG.

### 1.2 Extensions

**Mini-Batch.** Our result in this paper trivially extends to the mini-batch setting: if in each iteration we select $f_i(\cdot)$ for more than one random indices $i$, then we can accordingly define the gradient estimator and the result of this paper still holds. Note that the speed up that we obtain in this case comparing to gradient descent is $O((n/b)^{1/3})$ where $b$ is the mini-batch size. Therefore, the smaller $b$ is the better sequential running time we expect to see (which is also observed in our experiments).

**Other Smoothness Assumptions.** Our result generalizes to the setting when each $f_i(\cdot)$ enjoys a different smoothness parameter. In this setting one needs to select a random index $i \in [n]$ with a non-uniform distribution in order to obtain a faster running time. Our result also generalizes to the upper-lower smoothness setting. Instead of requiring each $f_i(\cdot)$ to be $L$-smooth, one can assume it is $L$-upper smooth and $l$-lower smooth, a notation introduced by [4]; in such a case, faster results can also be obtained using our same proof techniques.

**Sum-of-Non-Convex Objectives.** Our analogous proof also applies to the sum-of-non-convex setting which is the same Problem (1.1) except $f(\cdot)$ is guaranteed to be $\sigma$-strongly convex. Our obtained running time is $\widetilde{O}(n + \sqrt{n}L/\sigma)$ for SVRG. This is faster than the previous running time on SVRG which is $\widetilde{O}(n + L^2/\sigma^2)$, however, it is not faster than using SVRG+Catalyst which gives $\widetilde{O}(n + n^{3/4}\sqrt{L}/\sqrt{\sigma})$, see discussion in [4]. We do not include the details about this proof because it does not outperform SVRG+Catalyst.

**Other Variance-Reduction Methods.** Our proof generalizes to *all* variance-reduction methods. However, for brevity we demonstrate it only for the SVRG algorithm.

### 1.3 Other Related Works

For convex objectives, finding stationary points (or equivalently the global minimum) for problem (1.1) has received lots of attentions across machine learning and optimization communities; many first-order methods [8, 16, 23, 26] as well as their accelerations [1, 3, 18, 27, 28] have been proposed in the past a few years. Even in the case when $f(\cdot)$ is convex but each $f_i(\cdot)$ is non-convex, the problem (1.1) can be solved easily [4, 11, 24].

**Algorithm 1** Simplified SVRG method in the non-convex setting

**Input:** $x^{\phi}$ a starting vector, $S$ number of epochs, $m$ number of iterations per epoch, $\eta$ step length.

$x_0^1 \leftarrow x^{\phi}$
**for** $s \leftarrow 1$ **to** $S$ **do**
  $\widetilde{\mu} \leftarrow \nabla f(x_0^s)$
  **for** $k \leftarrow 0$ **to** $m-1$ **do**
    Pick $i$ uniformly at random in $\{1, \cdots, n\}$.
    $\widetilde{\nabla} \leftarrow \nabla f_i(x_k^s) - \nabla f_i(x_0^s) + \widetilde{\mu}$
    $x_{k+1}^s = x_k^s - \eta\widetilde{\nabla}$
  **end for**
  $x_0^{s+1} \leftarrow x_m^s$
**end for**
**return** $x_k^s$ for some random $s \in \{1, 2, \ldots, S\}$ and random $k \in \{1, 2, \ldots, m\}$.

The results of Li and Lin [17] and Ghadimi and Lan [13] unify the theory of non-convex and convex optimization in the following sense. They provide general first-order schemes such that, if the parameters are tuned properly, the schemes can converge (1) as fast as gradient descent in terms of finding an approximate stationary point; and (2) as fast as accelerated gradient descent [20] in terms of minimizing the objective if the function is convex. For the class of (1.1), their methods are only as slow as GD; in contrast, in this paper we prove theoretical convergence that is strictly faster than GD, which is both interesting and unknown.

A few days after the first version of this paper appeared on arXiv, we became aware of another group of authors that have independently obtained essentially the same result [21, 22]. [4]

## 2 Notations and Algorithm

A differentiable function $f\colon \mathbb{R}^n \to \mathbb{R}$ is $L$-smooth if for all pairs $x, y \in \mathbb{R}^n$ it satisfies $\|\nabla f(x) - \nabla f(y)\| \le L\|x-y\|$. An equivalent definition says for for all $x, y \in \mathbb{R}^n$:

$$-\frac{L}{2}\|x-y\|^2 \le f(x) - f(y) - \langle \nabla f(y), x-y \rangle \le \frac{L}{2}\|x-y\|^2 \ . \tag{2.1}$$

The main body of this paper proves our result based on three *false* simplification assumptions 4.2, 4.3 and 4.4 for the sake of sketching the high-level intuitions and highlighting the differences between our proof and known results. Our formal convergence proof is quite technical and included in the full paper.

In this high-level proof, we consider Algorithm 1, a simplified version of our final SVRG method for the non-convex setting. Note that both the snapshot point and the starting iterate $x_0^s$ of the $s$-th epoch have been chosen as the last iterate of the previous epoch in Algorithm 1.

**Remark 2.1.** In our final proof, we instead choose $x_0^s$ to be a weighted average of the last $m^{2/3}$ iterates from the previous epoch. See Algorithm 2 in the full paper. We demonstrate in Section 6 that this also a better choice in practice.

[4] These results also address gradient dominated functions, for which our main theorem also applies as follows. A non-convex function $f(\cdot)$ is $\tau$-gradient dominated if $f(x) - f(x^*) \le \tau\|\nabla f(x)\|^2$ for every $x$. Since our main theorem implies one can obtain $x$ satisfying $\|\nabla f(x)\|^2 \le \frac{1}{2\tau}(f(x_0) - f(x^*))$ using $O\big(n + n^{2/3}L\tau\big)$ stochastic gradients, by repeating it $\log_2(1/\varepsilon)$ times, we obtain an $\varepsilon$-minimizer of $f(\cdot)$ in $O\big((n + n^{2/3}L\tau)\log(1/\varepsilon)\big)$ stochastic gradients.

Throughout this paper we denote by $x_k^s$ the $k$-th iterate of the $s$-th epoch, by $\nabla_k^s = \nabla f(x_k^s)$ the full gradient at this iterate, and by $\widetilde{\nabla}_k^s = \nabla f(x_0^s) + \nabla_i f(x_k^s) - \nabla_i f(x_0^s)$ the gradient estimator which clearly satisfies $\mathbb{E}_i[\widetilde{\nabla}_k^s] = \nabla_k^s$. We denote by $i_k^s$ the random index $i$ chosen at iteration $k$ of epoch $s$. We also denote by $(\sigma_k^s)^2 \stackrel{\text{def}}{=} \|\nabla_k^s - \widetilde{\nabla}_k^s\|^2$ the variance of the gradient estimator $\widetilde{\nabla}_k^s$. Under these notations, our simplified SVRG algorithm in Algorithm 1 simply performs update $x_{k+1}^s \leftarrow x_k^s - \eta\widetilde{\nabla}_k^s$ for a fixed step length $\eta > 0$ that shall be specified later.

Since we focus mostly on analyzing a single epoch, when it is clear from the context, we drop the superscript $s$ and denote by $x_k$, $i_k$, $\nabla_k$, $\widetilde{\nabla}_k$, $\sigma_k^2$ respectively for $x_k^s$, $i_k^s$, $\nabla_k^s$, $\widetilde{\nabla}_k^s$, $(\sigma_k^s)^2$. We also denote by $\blacktriangledown_k^2 \stackrel{\text{def}}{=} \|\nabla_k\|_2^2$, $\sigma_{i:j}^2 \stackrel{\text{def}}{=} \sum_{k=i}^j \sigma_k^2$ and $\blacktriangledown_{i:j}^2 \stackrel{\text{def}}{=} \sum_{k=i}^j \blacktriangledown_k^2$ for notational simplicity.

## 3 Two Useful Lemmas

We first observe two simple lemmas. The first one describes the expected objective decrease between two consecutive iterations. This is a standard step that is used in analyzing gradient descent for smooth functions, and we additionally take into account the variance of the gradient estimator.

**Lemma 3.1** (gradient descent). *If $x_{k+1} = x_k - \eta\widetilde{\nabla}_k$ for some gradient estimator $\widetilde{\nabla}_k$ satisfying $\mathbb{E}[\widetilde{\nabla}_k] = \nabla_k = \nabla f(x_k)$, and if the step length $\eta \le \frac{1}{L}$, we have*

$$f(x_k) - \mathbb{E}[f(x_{k+1})] \ge \frac{\eta}{2}\nabla_k^2 - \frac{\eta^2 L}{2}\mathbb{E}\big[\sigma_k^2\big] \ .$$

*Proof.* By the smoothness of the function, we have

$$\begin{aligned}\mathbb{E}[f(x_{k+1})] &\le f(x_k) + \mathbb{E}\big[\langle \nabla f(x_k), x_{k+1} - x_k\rangle\big] + \mathbb{E}\Big[\frac{L}{2}\|x_{k+1} - x_k\|^2\Big] \\ &= f(x_k) - \eta\|\nabla f(x_k)\|^2 + \frac{\eta^2 L}{2}\mathbb{E}\Big[\|\widetilde{\nabla}_k\|^2\Big] \\ &= f(x_k) - \eta\|\nabla f(x_k)\|^2 + \frac{\eta^2 L}{2}\mathbb{E}\Big[\|\nabla f(x_k)\|^2 + \|\widetilde{\nabla}_k - \nabla f(x_k)\|^2\Big] \ .\end{aligned}$$

This immediately yields Lemma 3.1 by using the assumption that $\eta \le \frac{1}{L}$. □

The next lemma follows from the classical analysis of mirror descent methods.[5] However, we make novel use of it on top of a non-convex but smooth function.

**Lemma 3.2** (mirror descent). *If $x_{k+1} = x_k - \eta\widetilde{\nabla}_k$ for some gradient estimator $\widetilde{\nabla}_k$ satisfying $\mathbb{E}[\widetilde{\nabla}_k] = \nabla_k = \nabla f(x_k)$, then for every $u \in \mathbb{R}^d$ it satisfies*

$$f(x_k) - f(u) \le \frac{\eta}{2}\big(\blacktriangledown_k^2 + \mathbb{E}[\sigma_k^2]\big) + \big(\frac{1}{2\eta} + \frac{L}{2}\big)\|x_k - u\|^2 - \frac{1}{2\eta}\mathbb{E}\big[\|x_{k+1} - u\|^2\big] \ .$$

*Proof.* We first write the following inequality which follows from classical mirror-descent analysis.

[5] Mirror descent is a terminology mostly used in optimization literature, see for instance the textbook [7]. In machine learning contexts, mirror-descent analysis is essentially identical to regret analysis. In our SVRG method, the descent step $x_{k+1}^s \leftarrow x_k^s - \eta\widetilde{\nabla}_k^s$ can be interpreted as a mirror descent step in the Euclidean space (see for instance [2]), and therefore mirror-descent analysis applies.

For every $u \in \mathbb{R}^d$, it satisfies

$$\begin{aligned}
&\langle \nabla_k, x_k - u\rangle = \mathbb{E}[\langle \widetilde{\nabla}_k, x_k - u\rangle] \\
&= \mathbb{E}[\langle \widetilde{\nabla}_k, x_k - x_{k+1}\rangle + \langle \widetilde{\nabla}_k, x_{k+1} - u\rangle] \\
&= \mathbb{E}[\langle \widetilde{\nabla}_k, x_k - x_{k+1}\rangle - \frac{1}{2\eta}\|x_k - x_{k+1}\|^2 + \frac{1}{2\eta}\|x_k - u\|^2 - \frac{1}{2\eta}\|x_{k+1} - u\|^2] \\
&\leq \mathbb{E}\big[\frac{\eta}{2}\|\widetilde{\nabla}_k\|^2 + \frac{1}{2\eta}\|x_k - u\|^2 - \frac{1}{2\eta}\|x_{k+1} - u\|^2\big] \ . \qquad (3.1)
\end{aligned}$$

Above, $\langle \widetilde{\nabla}_k, x_{k+1} - u\rangle = -\frac{1}{2\eta}\|x_k - x_{k+1}\|^2 + \frac{1}{2\eta}\|x_k - u\|^2 - \frac{1}{2\eta}\|x_{k+1} - u\|^2$ is known as the three-point equality of Bregman divergence (in the special case of Euclidean space). The only inequality is because we have $\frac{1}{2}\|a\|^2 + \frac{1}{2}\|b\|^2 \geq \langle a, b\rangle$ for any two vectors $a, b$.

Classically in convex optimization, one would lower bound the left hand side of (3.1) by $f(x_k) - f(u)$ using the convexity of $f(\cdot)$. We take a different path here because our objective $f(\cdot)$ is not convex. Instead, we use the lower smoothness property of function $f$ which is the first inequality in (2.1) to deduce that $\langle \nabla_k, x_k - u\rangle \geq f(x_k) - f(u) - \frac{L}{2}\|x_k - u\|^2$. Combining this with inequality (3.1), and taking into account $\mathbb{E}[\|\widetilde{\nabla}_k\|^2] = \|\nabla_k\|^2 + \mathbb{E}[\sigma_k^2]$ by the definition of variance, we finish the proof of Lemma 3.2. □

Our main theorem is motivated by the linear-coupling framework [2]. In particular, we linearly couple the above gradient and mirror descent lemmas, together with a variance upper-bound lemma described in the next section.

# 4 Upper Bounding the Variance

**High-Level Ideas.** The key idea behind all variance-reduction literatures (such as SVRG [16], SAGA [8], and SAG [23]) is to prove that the variance $\mathbb{E}[(\sigma_k^s)^2]$ decreases as $s$ or $k$ increases. However, the only known technique to achieve so is to upper bound $\mathbb{E}[(\sigma_k^s)^2]$ "essentially" by $O\big(f(x_k^s) - f(x^*)\big)$, the objective distance to the minimum. Perhaps the only exception is the work on sum-of-non-convex but strongly-convex objectives [4, 24], where the authors upper bound $\mathbb{E}[(\sigma_k^s)^2]$ by $O\big(\|x_k^s - x^*\|^2\big)$, the squared vector distance to the minimum. Such techniques fail to apply in our non-convex setting, because gradient-descent based methods do not necessarily converge to the global minimum.

We take a different path in this paper. We upper bound $\mathbb{E}[(\sigma_k^s)^2]$ by $O\big(\|x_k^s - x_0^s\|^2\big)$, the squared vector distance between the current vector $x_k^s$ and the first vector (i.e., the snapshot) $x_0^s$ of the current epoch $s$. This is certainly tighter than $O\big(\|x_k^s - x^*\|^2\big)$ from prior work.[6] Moreover, the less we move away from the snapshot, the better upper bound we obtain on the variance. This is conceptually different from all existing literatures.

Furthermore, we in turn argue that $\|x_k^s - x_0^s\|^2$ is at most some constant times $f(x_k^s) - f(x_0^s)$. To prove so, we wish to select $u = x_0^s$ in Lemma 3.2 and telescope it for multiple iterations $k$, ideally for all the iterations $k$ within the same epoch. This is possible for convex objectives but impossible for non-convex ones. More specifically, the ratio between $(1/2\eta + L/2)$ and $(1/2\eta)$ can be much larger than 1, preventing us from telescoping more than $O(1/\eta L)$ iterations in any meaningful manner (see (4.1)). In contrast, this ratio would be identical to 1 in the convex setting, or even smaller than 1 in the strongly convex setting. For this reason, we define $\eta = \frac{1}{m^{2/3}L}$, divide each epoch into $O(m^{1/3})$ subepochs each consisting of $O(m^{2/3})$ consecutive iterations. Now we can telescope

[6] This new technique has also been applied to convex settings recently [1].

Lemma 3.2 for all the iterations within a subepoch because $m^{2/3} \le O(1/\eta L)$. Finally, we use vector inequalities (see (4.5)) to combine these upper bounds for sub-epochs into an upper bound on the entire epoch.

**Technical Details.** We choose $\eta = \frac{1}{m_0 L}$ for some parameter $m_0$ that divides $m$. We will set $m_0$ to be $m^{2/3}$ and the reason will become clear at the end of this section. Define $d = m/m_0$ so an epoch is divided into $d$ sub-epochs.

We make the following parameter choices

**Definition 4.1.** *Define $\beta_0 = 1$ and $\beta_t \stackrel{\text{def}}{=} (1+\eta L)^{-t} = (1+1/m_0)^{-t}$ for $t = 1, \dots, m_0 - 1$. We have $1 \ge \beta_t \ge 1/e > 1/3$.*

For each $k = 0, 1, \dots, m - m_0$, we sum up Lemma 3.2 for iterations $k, k+1, \dots, k+m_0-1$ with multiplicative weights $\beta_0, \dots, \beta_{m_0-1}$ respectively. The norm square terms shall telescope in this summation, and we arrive at

$$\sum_{t=0}^{m_0-1} \beta_t\big(f(x_{k+t}) - f(u)\big) \le \frac{\eta}{2} \sum_{t=0}^{m_0-1} \beta_t\big(\blacktriangledown_{k+t}^2 + \sigma_{k+t}^2\big) + \big(\frac{1}{2\eta} + \frac{L}{2}\big)\|x_k - u\|^2 - \frac{\beta_{m_0-1}}{2\eta}\|x_{k+m_0} - u\|^2 \ . \tag{4.1}$$

**Simplification 4.2.** *Since the weights $\beta_0, \dots, \beta_{m_0-1}$ are within each other by a constant factor, let us assume for simplicity that they are all equal to 1.*

If we choose $u = x_k$ and assume $\beta_t = 1$ for all $t$, we can rewrite (4.1) as

$$\frac{1}{m_0}\sum_{t=0}^{m_0-1} \big(f(x_{k+t}) - f(x_k)\big) \le \frac{\eta}{2}\frac{1}{m_0}\big(\blacktriangledown_{k:k+m_0-1}^2 + \sigma_{k:k+m_0-1}^2\big) - \frac{1}{6\eta m_0}\|x_{k+m_0} - x_k\|^2 \ . \tag{4.2}$$

**Simplification 4.3.** *Since the left hand side of (4.2) is describing the average objective value $f(x_k), f(x_{k+1}), \dots, f(x_{k+m_0-1})$ which is hard to analyze, let us assume for simplicity that it can be replaced with the* last *iteration in this subepoch, that is*

$$f(x_{k+m_0}) - f(x_k) \le \frac{\eta}{2}\frac{1}{m_0}\big(\blacktriangledown_{k:k+m_0-1}^2 + \sigma_{k:k+m_0-1}^2\big) - \frac{1}{6\eta m_0}\|x_{k+m_0} - x_k\|^2 \ . \tag{4.3}$$

Using the above inequality we provide a novel upper bound on the variance of the gradient estimator:

$$\begin{aligned} \mathbb{E}_{i_t}\big[\sigma_t^2\big] &= \mathbb{E}_{i_t}\Big[\big\|\big(\nabla f_{i_t}(x_t) - \nabla f_{i_t}(x_0)\big) - \big(\nabla f(x_t) - \nabla f(x_0)\big)\big\|^2\Big] \\ &\le \mathbb{E}_{i_t}\Big[\big\|\nabla f_{i_t}(x_t) - \nabla f_{i_t}(x_0)\big\|^2\Big] \\ &\le L^2\|x_t - x_0\|^2 \ . \end{aligned} \tag{4.4}$$

Above, the first inequality is because for any random vector $\zeta \in \mathbb{R}^d$, it holds that $\mathbb{E}\|\zeta - \mathbb{E}\zeta\|^2 = \mathbb{E}\|\zeta\|^2 - \|\mathbb{E}\zeta\|^2$, and the second inequality is by the smoothness of each $f_i(\cdot)$.

In particular, for $t = m$, we can upper bound $\sigma_m^2$ using (4.4) and multiple times of (4.3):

$$\begin{aligned} \mathbb{E}[\sigma_m^2] \le L^2\mathbb{E}\big[\|x_m - x_0\|^2\big] &\le L^2 d\,\mathbb{E}\big[\|x_m - x_{m-m_0}\|^2 + \|x_{m-m_0} - x_{m-2m_0}\|^2 + \cdots + \|x_{m_0} - x_0\|^2\big] \\ &\le L^2 d\,\mathbb{E}\Big[6\eta m_0\big(f(x_0) - f(x_m)\big) + 3\eta^2\big(\blacktriangledown_{0:m-1}^2 + \sigma_{0:m-1}^2\big)\Big] \ . \end{aligned} \tag{4.5}$$

Above, the first inequality follows from the vector inequality $\|v_1 + \cdots + v_d\|^2 \le d\big(\|v_1\|^2 + \cdots + \|v_d\|^2\big)$, and the second inequality follows from (4.3).

**Simplification 4.4.** *Suppose that (4.5) holds not only for $\sigma_m^2$ but actually for all $\sigma_0^2, \dots, \sigma_{m-1}^2$, then it satisfies*

$$\frac{1}{m}\mathbb{E}[\sigma_{0:m-1}^2] \le L^2 d\mathbb{E}\Big[6\eta m_0\big(f(x_0) - f(x_m)\big) + 3\eta^2\big(\blacktriangledown_{0:m-1}^2 + \sigma_{0:m-1}^2\big)\Big] \ . \tag{4.6}$$

As long as $3\eta^2 L^2 d \le \frac{1}{2m}$, (4.6) further implies

$$\frac{1}{m}\mathbb{E}[\sigma_{0:m-1}^2] \le 12\eta m_0 L^2 d \cdot \mathbb{E}\big[f(x_0) - f(x_m) + \frac{\eta}{2m_0}\blacktriangledown_{0:m-1}^2\big] \ . \tag{4.7}$$

This concludes our goal in this section which is to provide an upper bound (4.7) on the (average) variance by a constant times the objective difference $f(x_0) - f(x_m)$.

# 5 Final Theorem

By applying the gradient descent guarantee Lemma 3.1 to the entire epoch. We obtain that

$$f(x_0) - \mathbb{E}[f(x_m)] \ge \mathbb{E}\Big[\frac{\eta}{2}\blacktriangledown_{0:m-1}^2 - \frac{\eta^2 L}{2}\sigma_{0:m-1}^2\Big] \ .$$

Combining this with the variance upper bound (4.7), we immediately have

$$f(x_0) - \mathbb{E}[f(x_m)] \ge \frac{\eta}{2}\mathbb{E}[\blacktriangledown_{0:m-1}^2] - 6\eta^3 m_0 m L^3 d \cdot \mathbb{E}[f(x_0) - f(x_m) + \frac{\eta}{2m_0}\blacktriangledown_{0:m-1}^2] \ . \tag{5.1}$$

In other words, as long as $6\eta^3 m_0 m L^3 d \le \frac{1}{2}$, we arrive at

$$f(x_0) - \mathbb{E}[f(x_m)] \ge \frac{\eta}{6}\mathbb{E}[\blacktriangledown_{0:m-1}^2] \ . \tag{5.2}$$

Note that (5.2) is only for a single epoch and can be written as $f(x_0^s) - \mathbb{E}[f(x_m^s)] \ge \frac{\eta}{6}\mathbb{E}[\sum_{t=0}^{m-1}\|\nabla f(x_t^s)\|^2]$ in the general notation. Therefore, we can telescope it over all epochs $s = 1, 2, \dots, S$. Since we have chosen $x_0^s$, the initial vector in epoch $s$, to be $x_m^{s-1}$, the last vector of the previous epoch, we obtain that

$$\frac{1}{Sm}\sum_{s=1}^{S}\sum_{t=0}^{m-1}\mathbb{E}\big[\|\nabla f(x_t^s)\|^2\big] \le \frac{6}{\eta Sm}(f(x_0^1) - f(x_m^S)) \le \frac{6(f(x^\phi) - \min_x f(x))}{\eta Sm} \ .$$

At this point, if we randomly select $s \in [S]$ and $t \in [m]$ at the end of the algorithm, we conclude that

$$\mathbb{E}[\|\nabla f(x_t^s)\|^2] \le \frac{6(f(x^\phi) - \min_x f(x))}{\eta Sm} \ .$$

(We remark here that selecting an average iterate to output is a common step also used by GD or SGD for non-convex optimization. This step is often unnecessarily in practice.)

Finally, let use choose the parameters properly. We simply let $m = n$ be the epoch length. Since we have required $3\eta^2 L^2 d \le \frac{1}{2m}$ and $6\eta^3 m_0 m L^3 d \le \frac{1}{2}$ in the previous section, and both these requirements can be satisfied when $m_0^3 \ge 12m^2$, we set $m_0 = \Theta(m^{2/3}) = \Theta(n^{2/3})$. Accordingly $\eta = \frac{1}{m_0 L} = \Theta\big(\frac{1}{n^{2/3}L}\big)$. In sum,

**Theorem 5.1.** *Under the simplification assumptions 4.2, 4.3 and 4.4, by choosing* $m = n$ *and* $\eta = \Theta\big(\frac{1}{n^{2/3}L}\big)$*, the produced output* $x$ *of Algorithm 1 satisfies that*[7]

$$\mathbb{E}[\|\nabla f(x)\|^2] \le O\Big(\frac{L(f(x^\phi) - \min_x f(x))}{Sn^{1/3}}\Big) \enspace .$$

*In other words, to obtain a point* $x$ *satisfying* $\|\nabla f(x)\|^2 \le \varepsilon$*, the total number of iterations needed for Algorithm 1 is*

$$Sn = O\Big(\frac{n^{2/3}L(f(x^\phi) - \min_x f(x))}{\varepsilon}\Big) \enspace .$$

The amortized per-iteration complexity of SVRG is at most twice of SGD. Therefore, this is a factor of $\Omega(n^{1/3})$ faster than the full gradient descent method on solving (1.1).

# 6 Experiments

## 6.1 Empirical Risk Minimization with Non-Convex Loss

We consider binary classification on four standard datasets that can be found on the LibSVM website [10]:

- the adult (a9a) dataset (32,561 training samples, 16,281 testing samples, and 123 features).
- the web (w8a) dataset (49,749 training samples, 14,951 testing samples, and 300 features).
- the rcv1 (rcv1.binary) dataset (20,242 training samples, 677,399 testing samples, and 47,236 features).
- the mnist (class 1) dataset (60,000 training samples, 10,000 testing samples, and 780 features,

**Accuracy Experiment.** In the first experiment we apply SVRG on training the $\ell_2$-regularized ERM problem with six loss functions: logistic loss, squared loss, smoothed hinge loss (with smoothing parameters 0.01, 0.1 and 1 resp.), and smoothed zero-one loss (also known as sigmoid loss).[8] We wish to see how non-convex loss compares to convex ones in terms of testing accuracy (and thus in terms of the generalization error).

For each of the four datasets, we also randomly flip 1/4 fraction, 1/8 fraction, or zero fraction of the training example labels. The purpose of this manipulation is to introduce outliers to the training set. We therefore have $4 \times 3 = 12$ datasets in total. We choose epoch length $m = 2n$ as suggested by the paper SVRG for ERM experiments, and use the simple Algorithm 1 for both convex and non-convex loss functions.

We present the accuracy results partially in Figure 1 (and the full version can be found in Figure 4 in the appendix). The $y$-axis represents the classification testing accuracy, and the $x$-axis represents the number of passes to the dataset. (Each iteration of SVRG counts as $1/n$ pass and the full-gradient computation of SVRG counts as 1 pass.)

These plots are produced based on a *fair and careful* parameter-tuning and parameter-validation procedure that can be described in the following four steps. Step I: for each of the 12 datasets, we partition the training samples randomly into a training set of size 4/5 and a validation set of size 1/5. Step II: for each of the 12 datasets and each loss function, we enumerate over 10 choices of $\lambda$, the regularization parameter. For each $\lambda$, we tune SVRG on the *training set* with different step lengths $\eta$ and choose the best $\eta$ that gives the fastest training speed. Step III: for each of the 12 datasets and each loss function, we tune the best $\lambda$ using the *validation set*. That is, we use the selected $\eta$ from Step II to train the linear predictor, and apply it on the validation set to obtain the testing accuracy. We then select the $\lambda$ that gives the best testing accuracy for the validation set. Step IV: for each of the 12 datasets and each loss function, we apply the validated linear predictor to the *testing set* and present it in Figure 1 and Figure 4.

[7] Like in SGD, one can easily apply a Markov inequality to conclude that with probability at least 2/3 we have the same asymptotic upper bound on the deterministic quantity $\|\nabla f(x)\|^2$.

[8] For the sigmoid loss, we scale it properly so that its smoothness parameter is exactly 1. Unlike hinge loss, it is unnecessary to consider sigmoid loss with different smoothness parameters: one can carefully verify that by scaling up or down the weight of the $\ell_2$ regularizer, it is equivalent to changing the smoothness of the sigmoid loss.

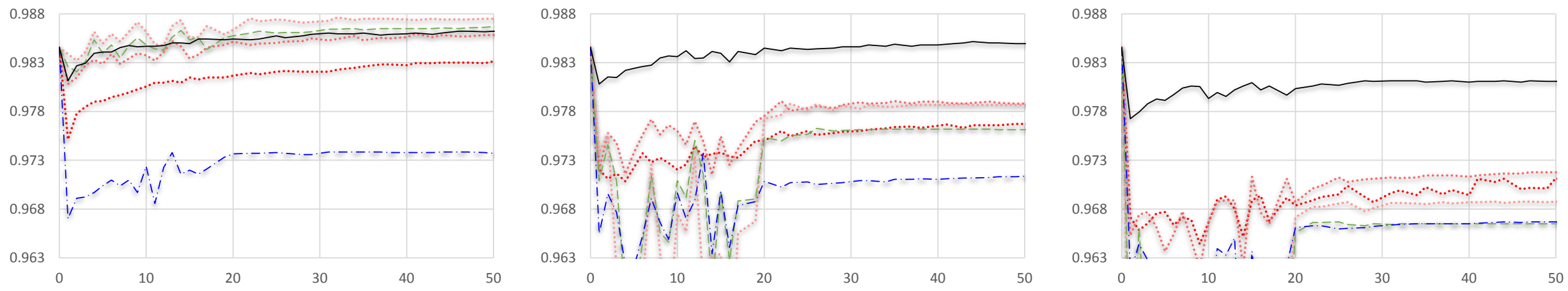

(a) web, no flip (b) web, flip 1/8 fraction (c) web, flip 1/4 fraction

Figure 1: Testing accuracy comparison on SVRG for $\ell_2$-regularized ERM with different loss functions. The full plots for all the 4 datasets can be found in Figure 4 in the appendix. **Black solid curves** represent sigmoid loss, blue dash curves represent square loss, green dash-dotted curves represent logistic loss, and the three red dotted curves represent hinge loss with 3 different smoothing parameters.

We make the following observations from this experiment.

- Although sigmoid loss is only comparable to hinge loss or logistic loss for "no flip" datasets, however, when the input has a lot of outliers (see "flip 1/8" and "flip 1/4"), sigmoid loss is undoubtedly the best choice. Square loss is almost always dominated because it is not necessarily a good choice for binary classification.

- The running time needed for SVRG on these datasets are quite comparable, regardless of the loss function being convex or not.

**Running-Time Experiment.** In this second experiment, we fix the regularization parameter $\lambda$ and compare the training objective convergence between SGD and SVRG for sigmoid loss only.[9] We choose four different $\lambda$ per dataset and present our plots partially in Figure 2 (and the full plots can be found in Figure 5 in the appendix). In these plots, the $y$-axis represents the training objective value, and the $x$-axis represents the number of passes to the dataset.

For a fair comparison we need to tune the step length $\eta$ for each dataset and each choice of $\lambda$. For SGD, we enumerate over polynomial learning rates $\eta_k = \alpha \cdot (1 + k/n)^\beta$ where $k$ is the number of iterations passed; we have made 10 choices of $\alpha$, considered $\beta = 0, 0.1, \dots, 1.0$, and selected the learning rate that gives a fastest convergence. For SVRG, we first consider the vanilla SVRG using a constant $\eta$ throughout all iterations, and select the best $\eta$ that gives the fastest convergence. This curve is presented in dashed blue in Figure 5. We also implement SVRG with polynomial learning rates $\eta_k = \alpha \cdot (1 + k/n)^\beta$ and tune the best $\alpha, \beta$ parameters and present the results in dashed black curves in Figure 5.

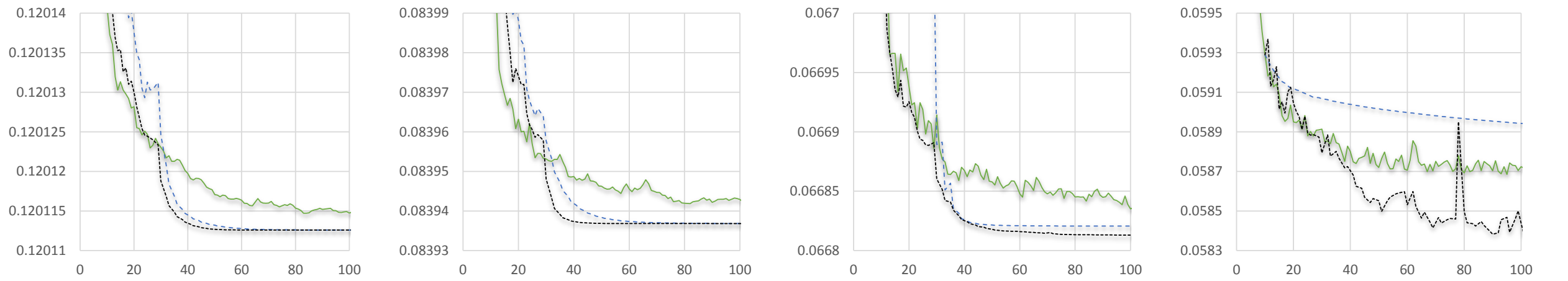

(a) web, $\lambda = 10^{-3}$ (b) web, $\lambda = 10^{-4}$ (c) web, $\lambda = 10^{-5}$ (d) web, $\lambda = 10^{-6}$

Figure 2: Training error comparison between SGD and SVRG on $\ell_2$-regularized ERM with sigmoid loss. The full plots can be found in Figure 5 in the appendix. The best-tuned SGD is presented in solid green, the best-tuned SVRG with constant step length is presented in dashed blue, and the best-tuned SVRG is presented in **doted black**.

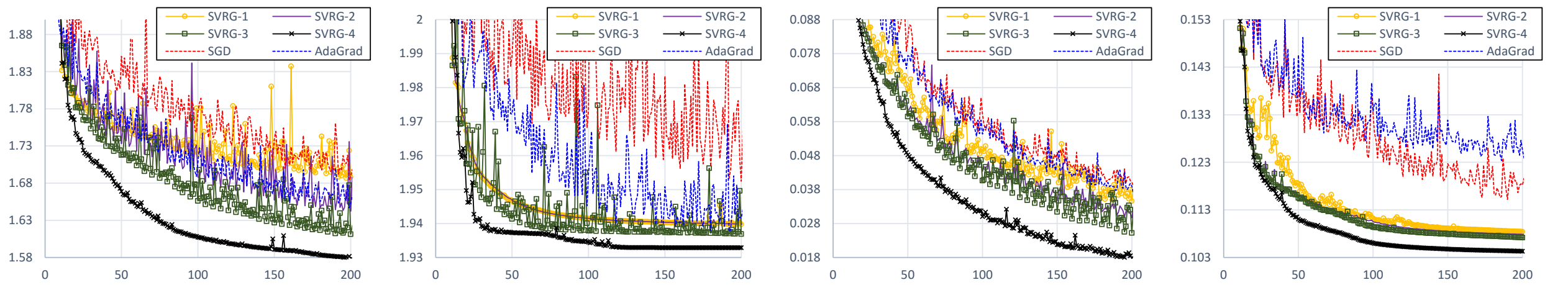

(a) cifar-10, $\lambda = 0$ (b) cifar-10, $\lambda = 10^{-3}$ (c) mnist, $\lambda = 0$ (d) mnist, $\lambda = 10^{-4}$

Figure 3: Training Error Comparison on neural nets.

We make the following observations from this experiment.

- Consistent with theory, SVRG is not necessarily a better choice than SGD for large training error $\varepsilon$. However, SVRG enjoys a very fast convergence especially for small $\varepsilon$.
- The smaller $\lambda$ is, the more "non-convex" the objective function becomes. We see that SGD performs more poorer than SVRG in these cases.[10]
- With only one exception (dataset web with $\lambda = 10^{-6}$), choosing a polynomial learning rate does not seem necessary in terms of improving the running time for training ERM problems with non-convex loss.
- Although not presented in Figure 5, the best-tuned polynomial learning rates for SVRG have smaller exponents $\beta$ as compared to SGD in each of the 16 plots.

### 6.2 Neural Network

We consider the multi-class (in fact, 10-class) classification problem on CIFAR-10 ($60,000$ training samples) and MNIST ($10,000$ training samples), two standard image datasets for neural net studies. We construct a toy two-layered neural net, with (1) 64 hidden nodes in the first layer, each connecting to a uniformly distributed 4x4 or 5x5 pixel block of the input image and having a smoothed relu (also known as softplus) activation function; (2) 10 output nodes on the second layer, fully connected to the first layer and each representing one of the ten classification outputs. We consider training such neural networks with the multi-class logistic loss that is a function on the 10 outputs and the correct label. For each of the two datasets, we consider both training the unregularized version, as well as the $\ell_2$ regularized version with weight $10^{-3}$ for CIFAR-10 and $10^{-4}$ for MNIST, two parameters suggested by [16].

[9]This experiment is the minimization problem with respect to all training samples since there is no need to perform validation.

[10]We note that the plots for different values $\lambda$ are presented with different vertical scales. For instance, at 100 passes of the dataset, the objective difference between SVRG and SGD is more than $2 \times 10^{-4}$ for $\lambda = 10^{-6}$ on dataset web, but less than $5 \times 10^{-6}$ for $\lambda = 10^{-3}$.

We implement two classical algorithms: stochastic gradient descent (SGD) with the best tuned polynomial learning rate and adaptive subGradient method (AdaGrad) of [9, 19] which is essentially SGD but with an adaptive learning rate. We choose a mini-batch size of 100 for both these methods. We consider four variants of SVRG, all of which use epoch length $m = 5n/b$ if $b$ is the mini-batch size:

- SVRG-1, the simple Algorithm 1 with a best tuned polynomial learning rate and $b = 100$.
- SVRG-2, our full Algorithm 2 with a best tuned polynomial learning rate and $b = 100$.[11]
- SVRG-3, using adaptive learning rate (similar to AdaGrad) on top of SVRG-2 with $b = 100$.
- SVRG-4, same as SVRG-3 but with $b = 16$.

Our training error performance is presented in Figure 3. We also include the testing accuracy in Figure 6 in the appendix. In these plots the $y$ axis represents the training objective value, and the $x$ axis represents the number of passes to the dataset. Each iteration of SGD or SVRG counts as $b/n$ pass of the dataset, and the snapshot full-gradient computation counts as 1 pass.[12] From the plots we clearly see a performance advantage for using SVRG-based algorithms as compared to SGD or AdaGrad. Furthermore, we observe that the following three features on top of SVRG could further improve its running time:

1. Comparing SVRG-2 with SVRG-1, we see that setting the epoch initial vector to be a weighted average of the last a few iterations of the previous epoch is recommended.

2. Comparing SVRG-3 with SVRG-2, we see that using adaptive learning rates comparing to tuning the best polynomial learning rate is recommended.

3. Comparing SVRG-4 with SVRG-3, we see that a smaller mini-batch size is recommended in terms of the total complexity. In contrast, reducing the mini-batch size is discouraged for SGD or AdaGrad because the variance could blow up and the performances would be decreased (this is also observed by our experiment but not included in the plots).

We hope that the above observations provide new insights for experimentalists working on deep learning.

[11] That is, we set the initial vector of each epoch to be weighted average of the last $(m/b)^{2/3}$ vectors from the previous epoch.

[12] The number of passes to the dataset is a traditional unit for comparing stochastic methods. For ERM problems, it is natural to count each iteration of SVRG as $b/n$ passes of the data rather than $2b/n$, because the computation of $\nabla f_i(\widetilde{x})$ is free if one efficiently stores $\nabla f_i(\widetilde{x})$ when the full gradient was computed at $\widetilde{x}$. However, after our paper has appeared online, we noticed this measurement may not be fair for SGD on training neural networks, because it is memory-inefficient to store $\nabla f_i(\widetilde{x})$ when $f_i$ comes from a large-scale neural network. For this reason, the *sequential* per-iteration cost of SVRG can be a factor $(2+1/5)/(1+1/5) = 11/6$ greater than SGD. Nevertheless, the extra cost on computing $\nabla f_i(\widetilde{x})$ is totally parallelizable (and can be viewed as doubling the mini-batch size), so this may not affect the GPU-based running time of SVRG by that much. We leave it a future work to run SVRG on large-scale network networks because it is beyond the scope of this paper.

**Algorithm 2** Our full SVRG method in the non-convex setting

**Input:** $x^\phi$ a starting vector, $S$ number of epochs, $m$ number of iterations per epoch, $m_0$ subepoch length, $\eta$ step length.

$x_0^1 \leftarrow x^\phi$
**for** $s \leftarrow 1$ **to** $S$ **do**
  $\widetilde{\mu} \leftarrow \nabla f(x_0^s)$
  **for** $k \leftarrow 0$ **to** $m-1$ **do**
    Pick $i$ uniformly at random in $\{1, \cdots, n\}$.
    $\widetilde{\nabla} \leftarrow \nabla f_i(x_k^s) - \nabla f_i(x_0^s) + \widetilde{\mu}$
    $x_{k+1}^s = x_k^s - \eta \widetilde{\nabla}$
  **end for**
  Define $\beta_0, \beta_1, \dots, \beta_{m_0-1}$ following Definition 4.1
  Select a random $m^s \in \{m, m-1, \dots, m-m_0+1\}$ with probability proportional to

$$\left\{\beta_{m_0-1}, \frac{10}{9}\beta_{m_0-1}, \frac{10}{9}(\beta_{m_0-1}+\beta_{m_0-2}), \dots, \frac{10}{9}(\beta_{m_0-1}+\cdots+\beta_1)\right\} .$$

  $x_0^{s+1} \leftarrow x_{m^s}^s$.
**end for**
**return** a vector uniformly at random from the set $\{x_{t-1}^s \,:\, s \in [S], t \in [m^s]\}$

## Acknowledgements

E. Hazan acknowledges support from the National Science Foundation grant IIS-1523815 and a Google research award. Z. Allen-Zhu acknowledges support from a Microsoft research award, no. 0518584.

# Appendix

## A Detailed Proof

In the detailed proof, we again first concentrating on analyzing a single epoch. We choose as before $\eta = \frac{1}{m_0 L}$ for some parameter $m_0$ that divides $m$. The natural choice of $m_0$ shall become clear at the end of this section, and would be set to around $m^{2/3}$. Define $d = m/m_0$ so an epoch is divided into $d$ sub-epochs. We denote by $x_{-m_0+1} = \cdots = x_{-1} \stackrel{\text{def}}{=} x_0$ for notational convenience, and similarly define $\nabla_{-m_0+1} = \cdots = \nabla_{-1} = \widetilde{\nabla}_{-m_0+1} = \cdots = \widetilde{\nabla}_{-1} = \sigma_{-m_0+1} = \cdots = \sigma_{-1} = 0$. Throughout this section, we also drop the expectation sign $\mathbb{E}[\cdot]$ for notational simplicity.

Our full algorithm for the non-convex setting is slightly different from our sketched proof, see Algorithm 2. Most importantly, instead of setting the first vector of each epoch to be the last iterate of the previous epoch, we set it to be a non-uniform random iterate in the last subepoch of the previous epoch. This step is crucial for our analysis to hold without the simplification assumptions.

### A.1 Upper Bounding the Variance

The following lemma is a simple counterpart to (4.4) in our sketched-proof section. It upper bounds the average variance inside an epoch by the average squared distances between vectors that are $m_0$ iterations away from each other.

**Lemma A.1.**

$$\sum_{t=0}^{m-1}\sigma_t^2 \le L^2 d^2 \sum_{t=0}^{m-1}\|x_{t+1}-x_{t+1-m_0}\|^2$$

*Proof.* Recall that we have that for every $t \in \{0,1,\dots,m-1\}$, we have

$$\sigma_t^2 \le L^2\|x_t - x_0\|^2 \le L^2 d(\|x_t - x_{t-m_0}\|^2 + \|x_{t-m_0}-x_{t-2m_0}\|^2 + \cdots)$$

Summing this up over all possible $t$'s, we have

$$\sum_{t=0}^{m-1}\sigma_t^2 \le L^2 d^2 \sum_{t=0}^{m-1}\|x_{t+1}-x_{t+1-m_0}\|^2 \enspace .$$

We emphasize that this analysis relies on our careful choice of $x_{-m_0+1} = \cdots = x_{-1} \stackrel{\mathrm{def}}{=} x_0$ which simplifies our notations. □

We next state a simple variant of Lemma 3.2 that allows negative indices:

**Lemma A.2.** *For every $k \in \{-m_0+1,\dots,m-m_0\}$, every $t\in\{0,\dots,m_0-1\}$, and every $u\in\mathbb{R}^d$, we have*

$$f(x_{k+t}) - f(u) \le \frac{\eta}{2}\big(\blacktriangledown_{k+t}^2 + \sigma_{k+t}^2\big) + \big(\frac{1}{2\eta}+\frac{L}{2}\big)\|x_{k+t}-u\|^2 - \frac{1}{2\eta}\|x_{k+t+1}-u\|^2$$

We define the same sequence of $\beta_0,\beta_1,\dots,\beta_{m_0-1}$ as before:

**Definition A.3.** *Define $\beta_0 = 1$ and $\beta_t \stackrel{\mathrm{def}}{=} (1+\eta L)^{-t} = (1+1/m_0)^{-t}$ for $t=1,\dots,m_0-1$. We have $1 \ge \beta_t \ge 1/e > 1/3$.*

By summing up Lemma A.2 with multiplicative ratios $\beta_t$ for each $t=0,1,\dots,m_0-1$, we arrive at the following lemma which is a counterpart of (4.1) in the sketched-proof section.

**Lemma A.4.** *For every $k\in\{-m_0+1,\dots,m-m_0\}$, and every $u$,*

$$\sum_{t=0}^{m_0-1}\beta_t\big(f(x_{k+t})-f(u)\big) \le \frac{\eta}{2}\sum_{t=0}^{m_0-1}\beta_t\big(\blacktriangledown_{k+t}^2+\sigma_{k+t}^2\big) + \big(\frac{1}{2\eta}+\frac{L}{2}\big)\|x_k-u\|^2 - \frac{\beta_{m_0-1}}{2\eta}\|x_{k+m_0}-u\|^2$$

*In particular, if we select $u = x_k$, we obtain that*

$$\sum_{t=1}^{m_0-1}\beta_t\big(f(x_{k+t})-f(x_k)\big) \le \frac{\eta}{2}\sum_{t=0}^{m_0-1}\beta_t\big(\blacktriangledown_{k+t}^2+\sigma_{k+t}^2\big) - \frac{1}{6\eta}\|x_{k+m_0}-x_k\|^2 \enspace . \tag{A.1}$$

The next lemma sums up (A.1) over all possible values of $k$. It can be viewed as a weighted, more sophisticated version of (4.6) in our sketched-proof section.

**Lemma A.5.** *As long as $m_0 \le \frac{1}{6\eta^2 L^2 d^2}$, we have*

$$\begin{aligned}
&\beta_{m_0-1}\big(f(x_{m-1}) - f(x_0) - 2\eta\blacktriangledown_{0:m-2}^2\big)\\
&+ (\beta_{m_0-1}+\beta_{m_0-2})\big(f(x_{m-2}) - f(x_0) - 2\eta\blacktriangledown_{0:m-3}^2\big) + \cdots\\
&+ (\beta_1+\cdots+\beta_{m_0-1})\big(f(x_{m-m_0+1}) - f(x_0) - 2\eta\blacktriangledown_{0:m-m_0}^2\big)\\
\le\ &\frac{\eta}{2}\beta_{m_0-1}\blacktriangledown_{m-1}^2 - \frac{1}{12\eta L^2 d^2}\sum_{t=0}^{m-1}\sigma_t^2 \enspace .
\end{aligned}$$

*Proof.* By carefully summing up (A.1) for $k \in \{-m_0+1, \dots, m-m_0\}$, we have that

$$
\begin{aligned}
&\beta_{m_0-1} f(x_{m-1}) + (\beta_{m_0-1}+\beta_{m_0-2}) f(x_{m-2}) + \cdots + (\beta_1 + \cdots + \beta_{m_0-1}) f(x_{m-m_0+1}) \\
&- (\beta_1 + \cdots + \beta_{m_0-1}) f(x_{-m_0+1}) - (\beta_2 + \cdots + \beta_{m_0-1}) f(x_{-m_0+2}) - \cdots - \beta_{m_0-1} f(x_{-1}) \\
&\le \frac{\eta}{2}\Big(\sum_{t=0}^{m_0-1}\beta_t\Big)\Big(\sum_{t=0}^{m-1}\sigma_t^2\Big) + \frac{\eta}{2}\Big(\sum_{t=0}^{m_0-1}\beta_t\Big)\Big(\sum_{t=0}^{m-m_0}\blacktriangledown_t^2\Big) \\
&+ \frac{\eta}{2}\Big(\beta_{m_0-1}\blacktriangledown_{m-1}^2 + (\beta_{m_0-1}+\beta_{m_0-2})\blacktriangledown_{m-2}^2 + \cdots + (\beta_1+\cdots+\beta_{m_0-1})\blacktriangledown_{m-m_0+1}^2\Big) \\
&- \frac{1}{6\eta}\sum_{t=0}^{m-1}\|x_{t+1} - x_{t+1-m_0}\|^2 \enspace .
\end{aligned}
$$

Using the fact that $x_{-m_0+1} = \cdots = x_{-1} \stackrel{\text{def}}{=} x_0$, we can rewrite the left hand side and get

$$
\begin{aligned}
&\beta_{m_0-1}\big(f(x_{m-1}) - f(x_0)\big) + (\beta_{m_0-1}+\beta_{m_0-2})\big(f(x_{m-2}) - f(x_0)\big) + \cdots \\
&+ (\beta_1 + \cdots + \beta_{m_0-1})\big(f(x_{m-m_0+1}) - f(x_0)\big) \\
&\le \frac{\eta}{2}\Big(\sum_{t=0}^{m_0-1}\beta_t\Big)\Big(\sum_{t=0}^{m-1}\sigma_t^2\Big) + \frac{\eta}{2}\Big(\sum_{t=0}^{m_0-1}\beta_t\Big)\Big(\sum_{t=0}^{m-m_0}\blacktriangledown_t^2\Big) \\
&+ \frac{\eta}{2}\Big(\beta_{m_0-1}\blacktriangledown_{m-1}^2 + \underbrace{(\beta_{m_0-1}+\beta_{m_0-2})\blacktriangledown_{m-2}^2 + \cdots + (\beta_1+\cdots+\beta_{m_0-1})\blacktriangledown_{m-m_0+1}^2}_{①}\Big) \\
&- \frac{1}{6\eta}\sum_{t=0}^{m-1}\|x_{t+1} - x_{t+1-m_0}\|^2 \enspace .
\end{aligned}
$$

Using the specific values of $\beta_t$'s, we can relax the terms in ① above and rewrite the above inequality as

$$
\begin{aligned}
&\beta_{m_0-1}\big(f(x_{m-1}) - f(x_0) - 2\eta\blacktriangledown_{m-2}^2\big) + (\beta_{m_0-1}+\beta_{m_0-2})\big(f(x_{m-2}) - f(x_0) - 2\eta\blacktriangledown_{m-3}^2\big) + \cdots \\
&+ (\beta_1 + \cdots + \beta_{m_0-1})\big(f(x_{m-m_0+1}) - f(x_0) - 2\blacktriangledown_{m-m_0}^2\big) \\
&\le \frac{\eta}{2}\Big(\sum_{t=0}^{m_0-1}\beta_t\Big)\Big(\sum_{t=0}^{m-1}\sigma_t^2\Big) + \underbrace{\frac{\eta}{2}\Big(\sum_{t=0}^{m_0-1}\beta_t\Big)\Big(\sum_{t=0}^{m-m_0-1}\blacktriangledown_t^2\Big)}_{②} \\
&+ \frac{\eta}{2}\beta_{m_0-1}\blacktriangledown_{m-1}^2 - \frac{1}{6\eta}\sum_{t=0}^{m-1}\|x_{t+1} - x_{t+1-m_0}\|^2 \enspace .
\end{aligned}
$$

Now we further relax the terms in ② above and further conclude that

$$
\begin{aligned}
&\beta_{m_0-1}\big(f(x_{m-1}) - f(x_0) - 2\eta\blacktriangledown_{0:m-2}^2\big) + (\beta_{m_0-1}+\beta_{m_0-2})\big(f(x_{m-2}) - f(x_0) - 2\eta\blacktriangledown_{0:m-3}^2\big) + \cdots \\
&+ (\beta_1 + \cdots + \beta_{m_0-1})\big(f(x_{m-m_0+1}) - f(x_0) - 2\eta\blacktriangledown_{0:m-m_0}^2\big) \\
&\le \frac{\eta}{2}\Big(\sum_{t=0}^{m_0-1}\beta_t\Big)\Big(\sum_{t=0}^{m-1}\sigma_t^2\Big) + \frac{\eta}{2}\beta_{m_0-1}\blacktriangledown_{m-1}^2 - \underbrace{\frac{1}{6\eta}\sum_{t=0}^{m-1}\|x_{t+1} - x_{t+1-m_0}\|^2}_{③} \enspace .
\end{aligned}
$$

Applying the variance bound Lemma A.1 on the summation ③, we have

$$\begin{aligned}&\beta_{m_0-1}\big(f(x_{m-1}) - f(x_0) - 2\eta \blacktriangledown^2_{0:m-2}\big) + (\beta_{m_0-1}+\beta_{m_0-2})\big(f(x_{m-2}) - f(x_0) - 2\eta \blacktriangledown^2_{0:m-3}\big) + \cdots \\ &+ (\beta_1 + \cdots + \beta_{m_0-1})\big(f(x_{m-m_0+1}) - f(x_0) - 2\eta \blacktriangledown^2_{0:m-m_0}\big) \\ \le\ & \frac{\eta}{2}\Big(\sum_{t=0}^{m_0-1}\beta_t\Big)\Big(\sum_{t=0}^{m-1}\sigma_t^2\Big) + \frac{\eta}{2}\beta_{m_0-1}\blacktriangledown^2_{m-1} - \frac{1}{6\eta L^2 d^2}\sum_{t=0}^{m-1}\sigma_t^2 \enspace.\end{aligned}$$

Finally, as long as $\sum_{t=0}^{m_0-1}\beta_t \le \frac{1}{6\eta^2 L^2 d^2}$ (which can be satisfied because $m_0 \le \frac{1}{6\eta^2 L^2 d^2}$), we have

$$\begin{aligned}&\beta_{m_0-1}\big(f(x_{m-1}) - f(x_0) - 2\eta \blacktriangledown^2_{0:m-2}\big) + (\beta_{m_0-1}+\beta_{m_0-2})\big(f(x_{m-2}) - f(x_0) - 2\eta \blacktriangledown^2_{0:m-3}\big) + \cdots \\ &+ (\beta_1 + \cdots + \beta_{m_0-1})\big(f(x_{m-m_0+1}) - f(x_0) - 2\eta \blacktriangledown^2_{0:m-m_0}\big) \\ \le\ & \frac{\eta}{2}\beta_{m_0-1}\blacktriangledown^2_{m-1} - \frac{1}{12\eta L^2 d^2}\sum_{t=0}^{m-1}\sigma_t^2 \enspace.\end{aligned}$$

This finishes the proof of Lemma A.5. □

## A.2 Objective Decrease using Gradient Descent

The following lemma is a variant of (5.1). However, instead of lower bounding the objective decrease $f(x_0)-f(x_m)$ for the entire epoch as in the sketched-proof section, we have to carefully lower bound a weighted sum of $f(x_0) - f(x_t)$ for $t \in \{m, m-1, \dots, m-m_0+1\}$, in order to make it consistent with the left hand side of Lemma A.5.

**Lemma A.6.**

$$\begin{aligned}&\beta_{m_0-1}\big(f(x_0) - f(x_m) - \frac{\eta}{2}\blacktriangledown^2_{0:m-1}\big) + \beta_{m_0-1}\big(f(x_0) - f(x_{m-1}) - \frac{\eta}{2}\blacktriangledown^2_{0:m-2}\big) \\ &+ (\beta_{m_0-1}+\beta_{m_0-2})\big(f(x_0) - f(x_{m-2}) - \frac{\eta}{2}\blacktriangledown^2_{0:m-3}\big) + \cdots \\ &+ (\beta_1 + \cdots + \beta_{m_0-1})\big(f(x_0) - f(x_{m-m_0+1}) - \frac{\eta}{2}\blacktriangledown^2_{0:m-m_0}\big) \\ \ge\ & -\frac{\eta^2 L m_0}{2}\sum_{t=0}^{m-1}\sigma_t^2\end{aligned}$$

*Proof.* For each $j = 1, 2, \dots, m_0$, by telescoping Lemma 3.1 across iterations $k = 0, 1, \dots, m-j$, we arrive at inequality

$$f(x_0) - f(x_{m-j+1}) \ge \frac{\eta}{2}\blacktriangledown^2_{0:m-j} - \frac{\eta^2 L}{2}\sum_{t=0}^{m-1}\sigma_t^2 \enspace.$$

Now we write down these $m_0$ inequalities separately, and multiply each of them by a positive weight:

$$\beta_{m_0-1} \times \Big\{ f(x_0) - f(x_m) \geq \frac{\eta}{2} \blacktriangledown^2_{0:m-1} - \frac{\eta^2 L}{2} \sum_{t=0}^{m-1} \sigma_t^2 \Big\}$$

$$\beta_{m_0-1} \times \Big\{ f(x_0) - f(x_{m-1}) \geq \frac{\eta}{2} \blacktriangledown^2_{0:m-2} - \frac{\eta^2 L}{2} \sum_{t=0}^{m-1} \sigma_t^2 \Big\}$$

$$(\beta_{m_0-2} + \beta_{m_0-1}) \times \Big\{ f(x_0) - f(x_{m-2}) \geq \frac{\eta}{2} \blacktriangledown^2_{0:m-3} - \frac{\eta^2 L}{2} \sum_{t=0}^{m-1} \sigma_t^2 \Big\} \quad \cdots$$

$$(\beta_1 + \cdots + \beta_{m_0-1}) \times \Big\{ f(x_0) - f(x_{m-m_0+1}) \geq \frac{\eta}{2} \blacktriangledown^2_{0:m-m_0} - \frac{\eta^2 L}{2} \sum_{t=0}^{m-1} \sigma_t^2 \Big\}$$

Summing these inequalities up, we obtain our desired inequality

$$\begin{aligned}
&\beta_{m_0-1}\big(f(x_0) - f(x_m) - \frac{\eta}{2} \blacktriangledown^2_{0:m-1}\big) + \beta_{m_0-1}\big(f(x_0) - f(x_{m-1}) - \frac{\eta}{2} \blacktriangledown^2_{0:m-2}\big) \\
&+ (\beta_{m_0-1} + \beta_{m_0-2})\big(f(x_0) - f(x_{m-2}) - \frac{\eta}{2} \blacktriangledown^2_{0:m-3}\big) + \cdots \\
&+ (\beta_1 + \cdots + \beta_{m_0-1})\big(f(x_0) - f(x_{m-m_0+1}) - \frac{\eta}{2} \blacktriangledown^2_{0:m-m_0}\big) \\
\geq& -\frac{\eta^2 L m_0^2}{2} \sum_{t=0}^{m-1} \sigma_t^2 \ .
\end{aligned}$$

□

## A.3 Final Theorem

Let us now put together Lemma A.5 and Lemma A.6, and derive the following lemma:

**Lemma A.7.** *As long as $6\eta^3 L^3 m_0^2 d^2 = 1/9$, we have*

$$\begin{aligned}
&\beta_{m_0-1}\big(f(x_0) - f(x_m) - \frac{\eta}{4} \blacktriangledown^2_{0:m-1}\big) + \frac{10\beta_{m_0-1}}{9}\big(f(x_0) - f(x_{m-1}) - \frac{\eta}{4} \blacktriangledown^2_{0:m-2}\big) \\
&+ \frac{10}{9}(\beta_{m_0-1} + \beta_{m_0-2})\big(f(x_0) - f(x_{m-2}) - \frac{\eta}{4} \blacktriangledown^2_{0:m-3}\big) + \cdots \\
&+ \frac{10}{9}(\beta_1 + \cdots + \beta_{m_0-1})\big(f(x_0) - f(x_{m-m_0+1}) - \frac{\eta}{4} \blacktriangledown^2_{0:m-m_0}\big) \geq 0 \ .
\end{aligned}$$

*Proof.* By directly combining Lemma A.5 and Lemma A.6, we have

$$\begin{aligned}
&\beta_{m_0-1}\big(f(x_0) - f(x_m) - \frac{\eta}{2} \blacktriangledown^2_{0:m-1}\big) + \beta_{m_0-1}\big(f(x_0) - f(x_{m-1}) - \frac{\eta}{2} \blacktriangledown^2_{0:m-2}\big) \\
&+ (\beta_{m_0-1} + \beta_{m_0-2})\big(f(x_0) - f(x_{m-2}) - \frac{\eta}{2} \blacktriangledown^2_{0:m-3}\big) + \cdots \\
&+ (\beta_1 + \cdots + \beta_{m_0-1})\big(f(x_0) - f(x_{m-m_0+1}) - \frac{\eta}{2} \blacktriangledown^2_{0:m-m_0}\big) \\
\geq& \ 12\eta L^2 d^2 \cdot \frac{\eta^2 L m_0^2}{2} \cdot \Big(\beta_{m_0-1}\big(f(x_{m-1}) - f(x_0) - 2\eta \blacktriangledown^2_{0:m-2}\big) \\
&+ (\beta_{m_0-1} + \beta_{m_0-2})\big(f(x_{m-2}) - f(x_0) - 2\eta \blacktriangledown^2_{0:m-3}\big) + \cdots \\
&+ (\beta_1 + \cdots + \beta_{m_0-1})\big(f(x_{m-m_0+1}) - f(x_0) - 2\eta \blacktriangledown^2_{0:m-m_0}\big) - \frac{\eta}{2}\beta_{m_0-1} \blacktriangledown^2_{m-1}\Big)
\end{aligned}$$

Suppose we have $12\eta L^2 d^2 \cdot \frac{\eta^2 L m_0^2}{2} = 6\eta^3 L^3 m_0^2 d^2 = 1/9$, then it satisfies that

$$\begin{aligned}
&\beta_{m_0-1}\big(f(x_0) - f(x_m) - \frac{\eta}{2}\blacktriangledown^2_{0:m-1}\big) + \beta_{m_0-1}\big(f(x_0) - f(x_{m-1}) - \frac{\eta}{2}\blacktriangledown^2_{0:m-2}\big) \\
&+ (\beta_{m_0-1} + \beta_{m_0-2})\big(f(x_0) - f(x_{m-2}) - \frac{\eta}{2}\blacktriangledown^2_{0:m-3}\big) + \cdots \\
&+ (\beta_1 + \cdots + \beta_{m_0-1})\big(f(x_0) - f(x_{m-m_0+1}) - \frac{\eta}{2}\blacktriangledown^2_{0:m-m_0}\big) \\
\geq\ & \frac{1}{9}\Big(\beta_{m_0-1}\big(f(x_{m-1}) - f(x_0) - 2\eta\blacktriangledown^2_{0:m-2}\big) \\
&+ (\beta_{m_0-1} + \beta_{m_0-2})\big(f(x_{m-2}) - f(x_0) - 2\eta\blacktriangledown^2_{0:m-3}\big) + \cdots \\
&+ (\beta_1 + \cdots + \beta_{m_0-1})\big(f(x_{m-m_0+1}) - f(x_0) - 2\eta\blacktriangledown^2_{0:m-m_0}\big) - \frac{\eta}{2}\beta_{m_0-1}\blacktriangledown^2_{m-1}\Big)
\end{aligned}$$

After rearranging, we have

$$\begin{aligned}
&\beta_{m_0-1}\big(f(x_0) - f(x_m) - \frac{\eta}{2}\blacktriangledown^2_{0:m-1}\big) + \frac{10\beta_{m_0-1}}{9}\big(f(x_0) - f(x_{m-1}) - \frac{\eta}{4}\blacktriangledown^2_{0:m-2}\big) \\
&+ \frac{10}{9}(\beta_{m_0-1} + \beta_{m_0-2})\big(f(x_0) - f(x_{m-2}) - \frac{\eta}{4}\blacktriangledown^2_{0:m-3}\big) + \cdots \\
&+ \frac{10}{9}(\beta_1 + \cdots + \beta_{m_0-1})\big(f(x_0) - f(x_{m-m_0+1}) - \frac{\eta}{4}\blacktriangledown^2_{0:m-m_0}\big) \geq -\frac{\eta}{18}\beta_{m_0-1}\blacktriangledown^2_{m-1} \ .
\end{aligned}$$

After relaxing the right hand side, we conclude that

$$\begin{aligned}
&\beta_{m_0-1}\big(f(x_0) - f(x_m) - \frac{\eta}{4}\blacktriangledown^2_{0:m-1}\big) + \frac{10\beta_{m_0-1}}{9}\big(f(x_0) - f(x_{m-1}) - \frac{\eta}{4}\blacktriangledown^2_{0:m-2}\big) \\
&+ \frac{10}{9}(\beta_{m_0-1} + \beta_{m_0-2})\big(f(x_0) - f(x_{m-2}) - \frac{\eta}{4}\blacktriangledown^2_{0:m-3}\big) + \cdots \\
&+ \frac{10}{9}(\beta_1 + \cdots + \beta_{m_0-1})\big(f(x_0) - f(x_{m-m_0+1}) - \frac{\eta}{4}\blacktriangledown^2_{0:m-m_0}\big) \geq 0 \ .
\end{aligned}$$

This finishes the proof of Lemma A.7. □

Lemma A.8 naturally implies that if we select a random stopping vector for this epoch, we have the following corollary which is a counterpart of (5.2) in our sketch-proof section:

**Corollary A.8.** *If we set $m^s$ to be a random variable in $\{m, m-1, \ldots, m-m_0+1\}$, with probabilities proportional to $\big\{\beta_{m_0-1}, \frac{10}{9}\beta_{m_0-1}, \frac{10}{9}(\beta_{m_0-1} + \beta_{m_0-2}), \ldots, \frac{10}{9}(\beta_{m_0-1} + \cdots + \beta_1)\big\}$, then Lemma A.7 implies that we have*

$$\mathbb{E}\big[f(x_0) - f(x_{m^s}) - \frac{\eta}{4}\blacktriangledown^2_{0:m^s-1}\big] \geq 0 \ .$$

Note that Corollary A.8 is only for a single epoch and can be written as

$$\mathbb{E}[f(x_0^s) - f(x_{m^s}^s)] \geq \frac{\eta}{4}\mathbb{E}\Big[\sum_{t=0}^{m^s-1} \|\nabla f(x_t^s)\|^2\Big]$$

in the general notation. Therefore, we are now ready to telescope it over all the epochs $s = 1, 2, \ldots, S$. Recall that we have chosen $x_0^s$, the initial vector in epoch $s$, to be $x_{m^{s-1}}^{s-1}$, the random stopping vector from the previous epoch. Therefore, we obtain that

$$\begin{aligned}
\frac{1}{m^1 + \cdots + m^S}\sum_{s=1}^{S}\sum_{t=0}^{m^s-1} \mathbb{E}\big[\|\nabla f(x_t^s)\|^2\big] &\leq \frac{4}{\eta S(m^1 + \cdots + m^S)}\big(f(x_0^1) - \mathbb{E}[f(x_{m^S}^S)]\big) \\
&\leq O\Big(\frac{f(x^\phi) - \min_x f(x)}{\eta S m}\Big) \ .
\end{aligned}$$

At this point, if we select uniformly at random an output vector $x$ from the set $\{x_{t-1}^s : s \in [S], t \in [m^s]\}$, we conclude that

$$\mathbb{E}[\|\nabla f(x_t^s)\|^2] \leq O\Big(\frac{f(x^\phi) - \min_x f(x)}{\eta S m}\Big) \ .$$

Finally, let use choose the parameters properly. We simply let $m = n$ be the epoch length. Since we have required $m_0 \leq \frac{1}{6\eta^2 L^2 d^2}$ and $6\eta^3 L^3 m_0^2 d^2 = 1/9$ in Lemma A.5 and Lemma A.7 respectively, and both these requirements can be satisfied when $m_0^3 \geq 54m^2$, we set $m_0 = \Theta(m^{2/3}) = \Theta(n^{2/3})$. Accordingly $\eta = \frac{1}{m_0 L} = O\big(\frac{1}{n^{2/3} L}\big)$. In sum,

**Theorem A.9** (Formal statement of Theorem 5.1)**.** *By choosing $m = n$ and $\theta = \Theta\big(\frac{1}{n^{2/3}L}\big)$, the produced output $x$ of Algorithm 2 satisfies that*[13]

$$\mathbb{E}[\|\nabla f(x)\|^2] \leq O\Big(\frac{L(f(x^\phi) - \min_x f(x))}{S n^{1/3}}\Big) \ .$$

*In other words, to obtain a point $x$ satisfying $\|\nabla f(x)\|^2 \leq \varepsilon$, the total number of iterations needed for Algorithm 1 is*

$$Sn = O\Big(\frac{n^{2/3} L(f(x^\phi) - \min_x f(x))}{\varepsilon}\Big) \ .$$

[13]Like in SGD, one can easily apply a Markov inequality to conclude that with probability at least 2/3 we have the same asymptotic upper bound on the deterministic quantity $\|\nabla f(x)\|^2$.

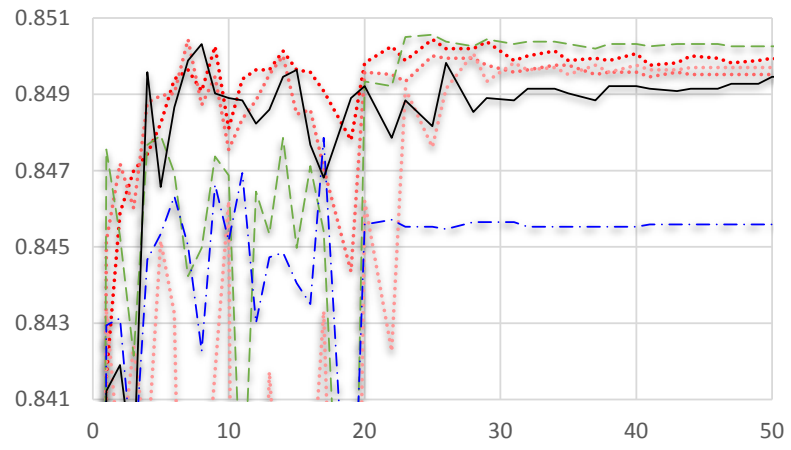

(a) adult, no flip

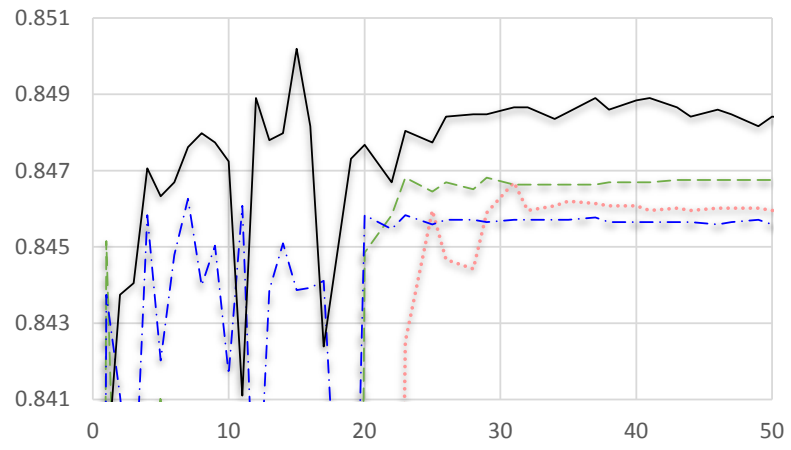

(b) adult, flip 1/8 fraction

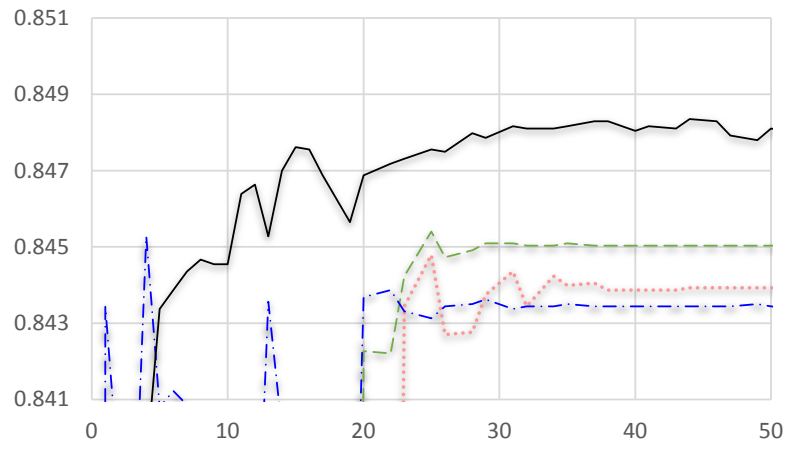

(c) adult, flip 1/4 fraction

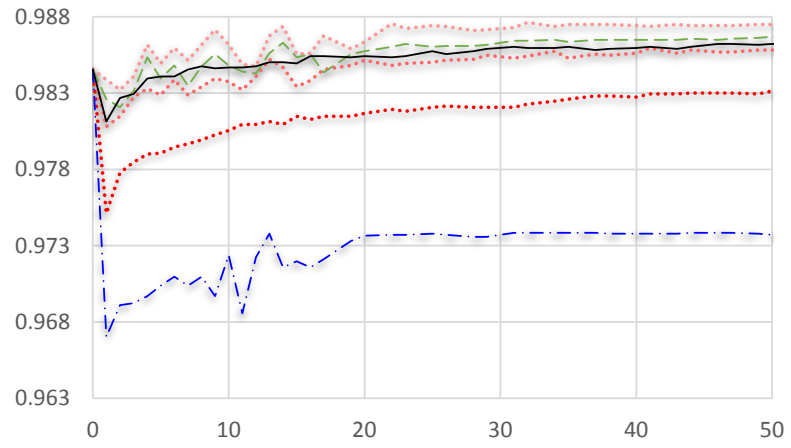

(d) web, no flip

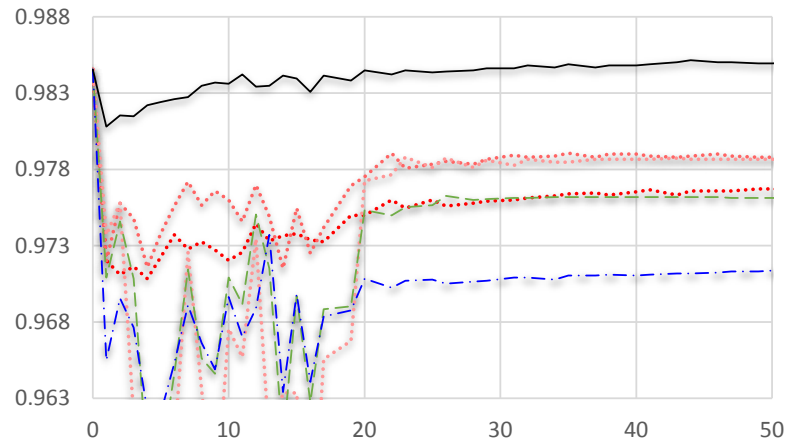

(e) web, flip 1/8 fraction

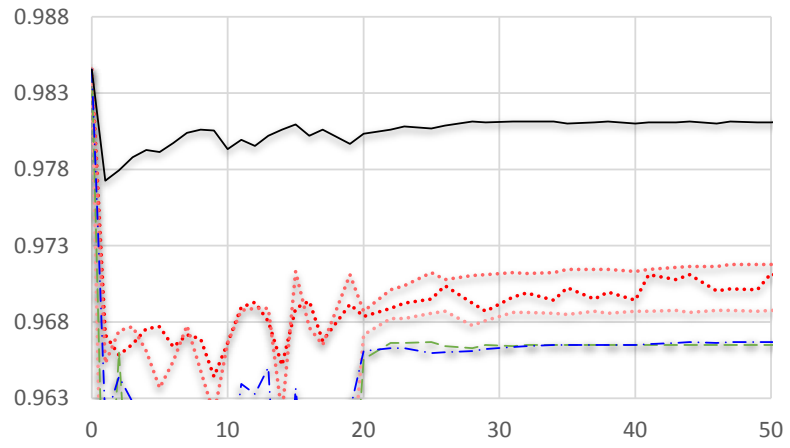

(f) web, flip 1/4 fraction

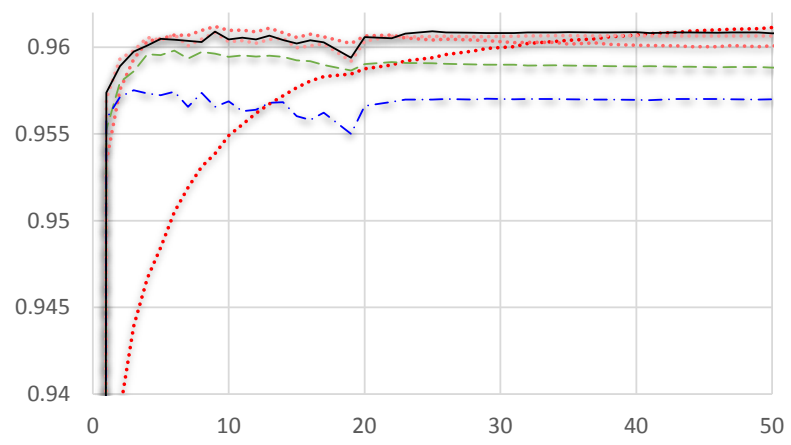

(g) rcv1, no flip

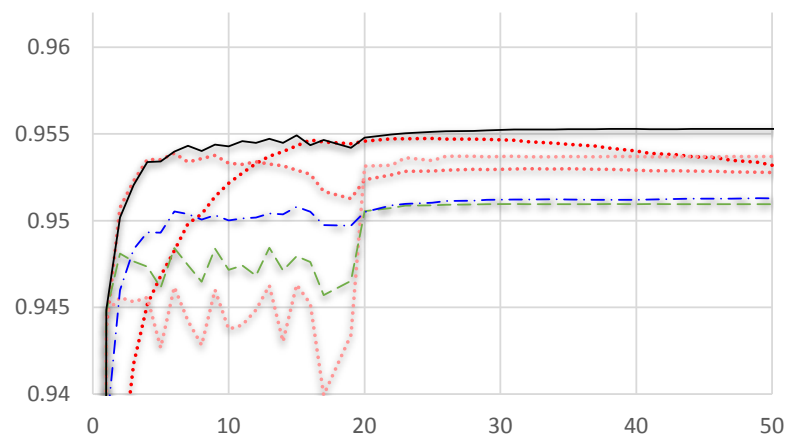

(h) rcv1, flip 1/8 fraction

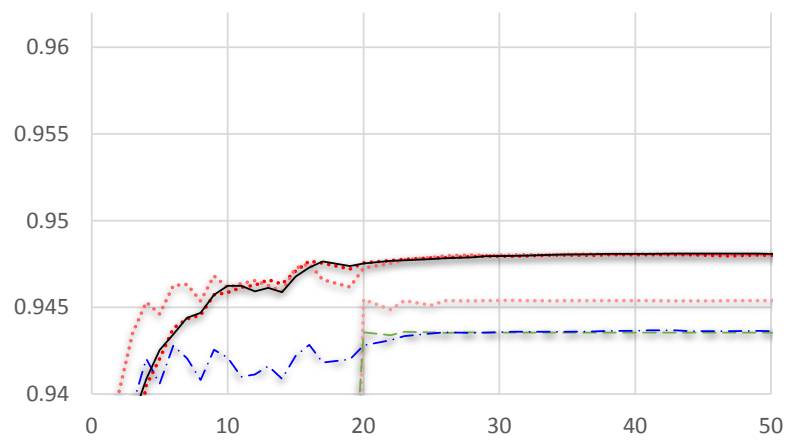

(i) rcv1, flip 1/4 fraction

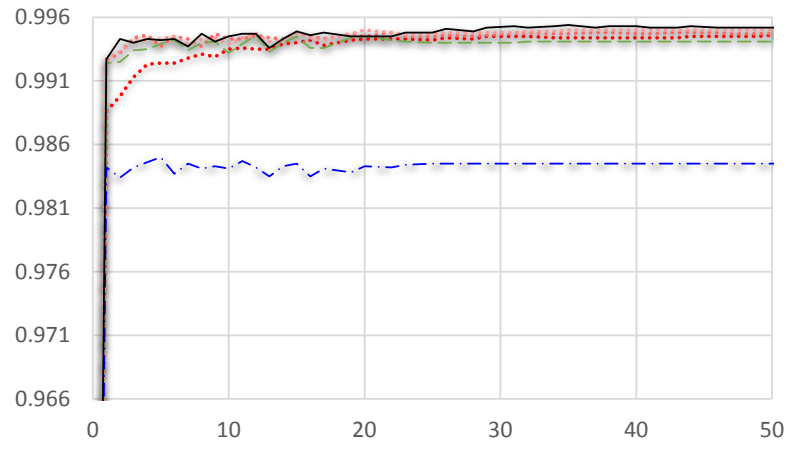

(j) mnist, no flip

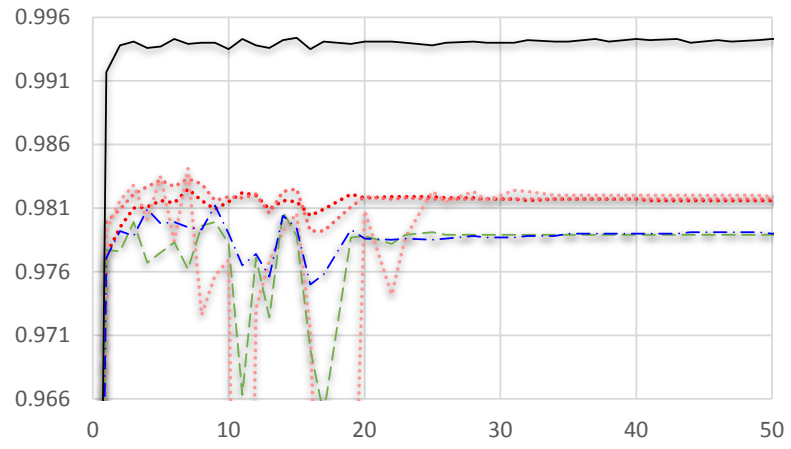

(k) mnist, flip 1/8 fraction

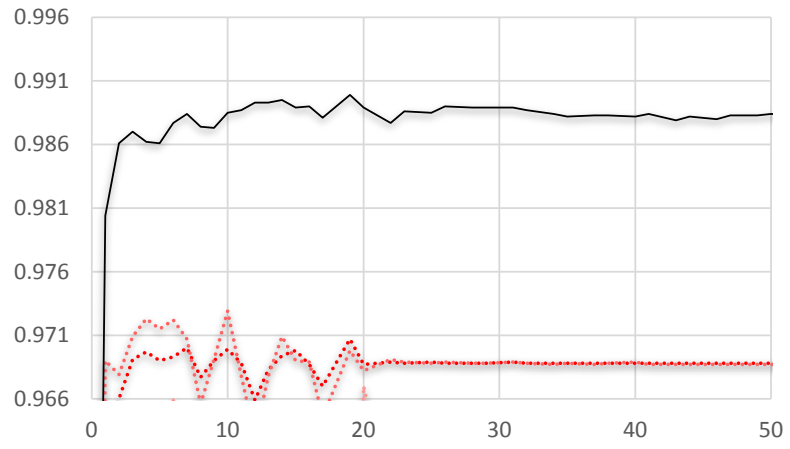

(l) mnist, flip 1/4 fraction

Figure 4: Testing accuracy comparison on SVRG for $\ell_2$-regularized ERM with different loss functions. **Black solid curves** represent sigmoid loss, blue dash curves represent square loss, green dash-dotted curves represent logistic loss, and the three red dotted curves represent hinge loss with 3 different smoothing parameters.

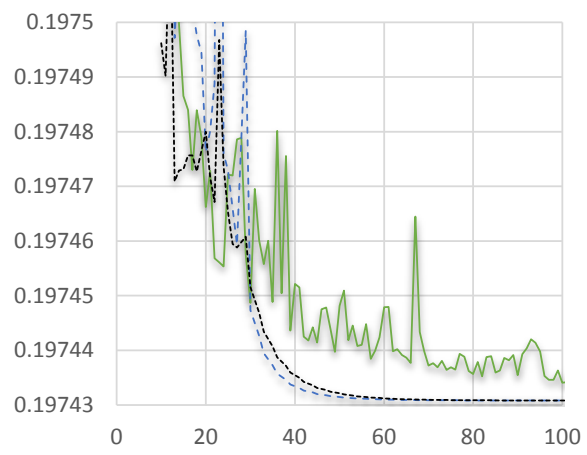

(a) adult, $\lambda = 10^{-3}$

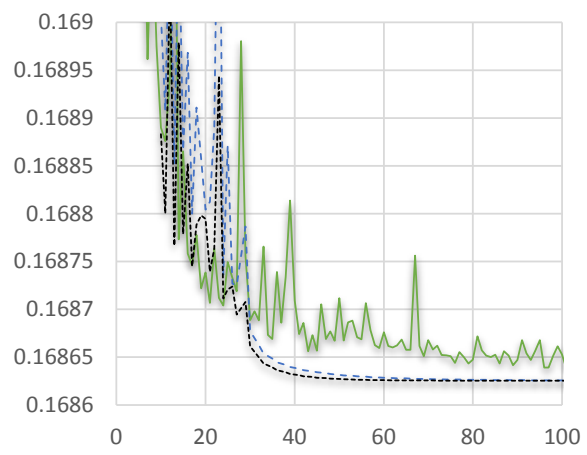

(b) adult, $\lambda = 10^{-4}$

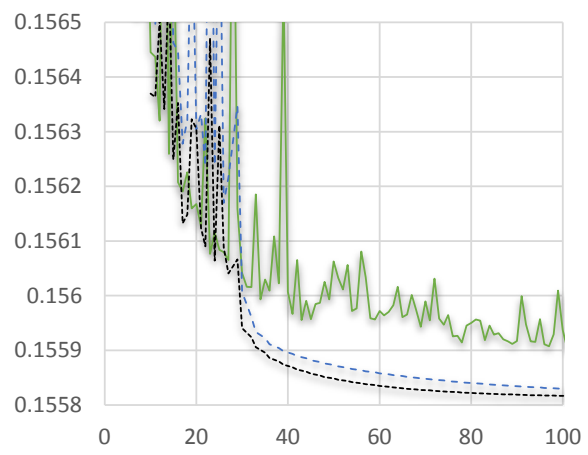

(c) adult, $\lambda = 10^{-5}$

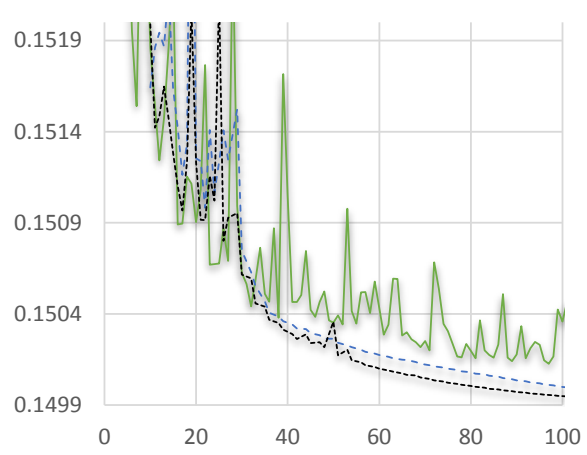

(d) adult, $\lambda = 10^{-6}$

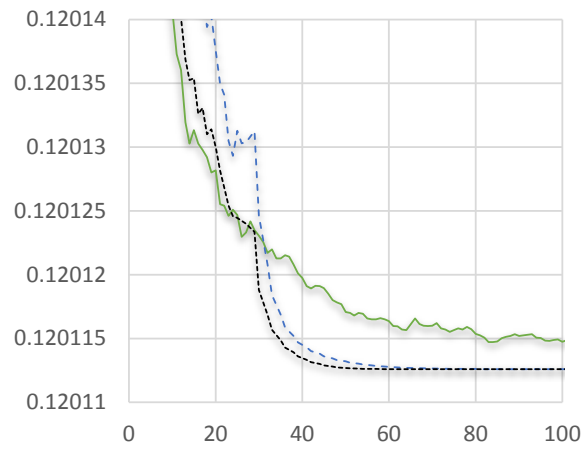

(e) web, $\lambda = 10^{-3}$

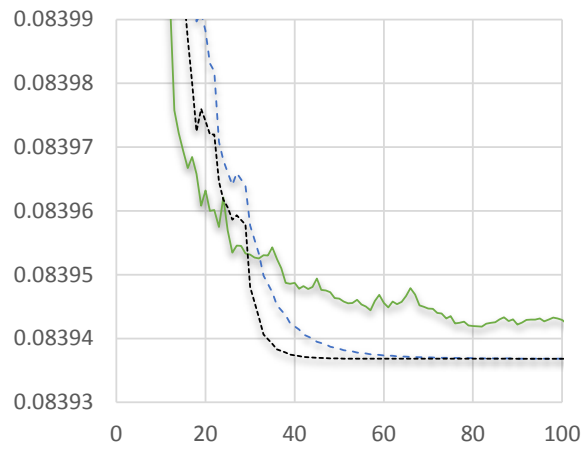

(f) web, $\lambda = 10^{-4}$

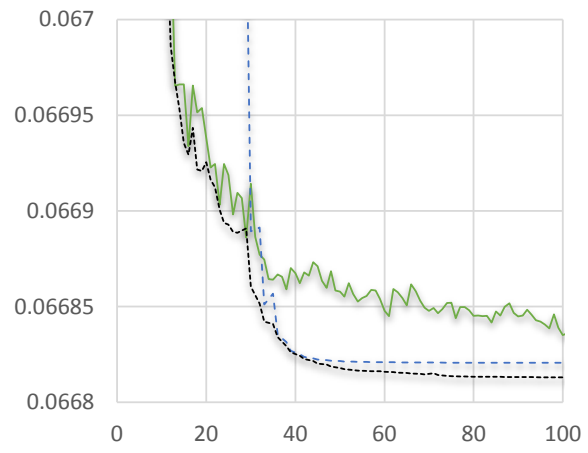

(g) web, $\lambda = 10^{-5}$

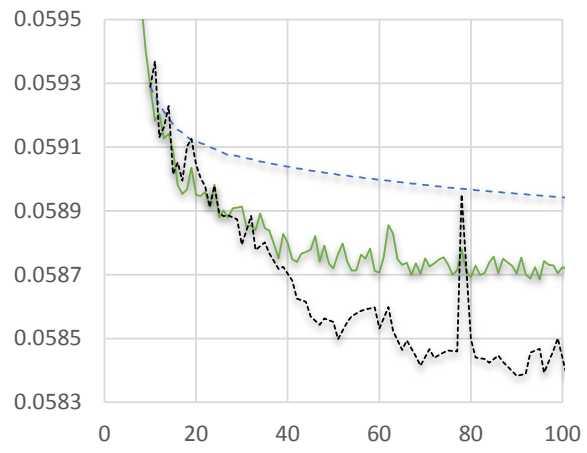

(h) web, $\lambda = 10^{-6}$

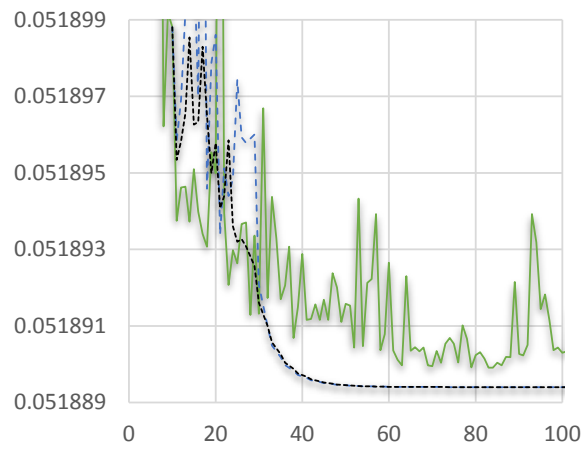

(i) mnist, $\lambda = 10^{-3}$

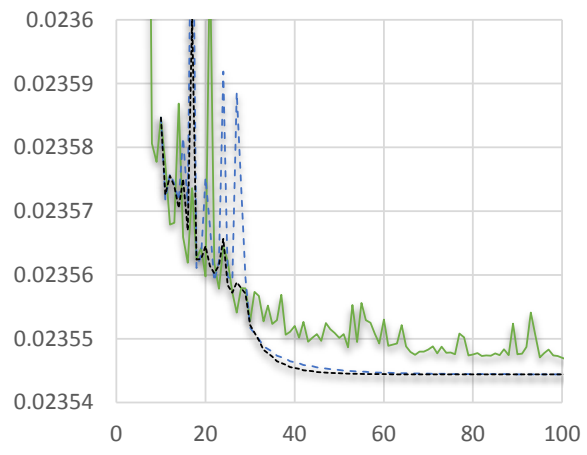

(j) mnist, $\lambda = 10^{-4}$

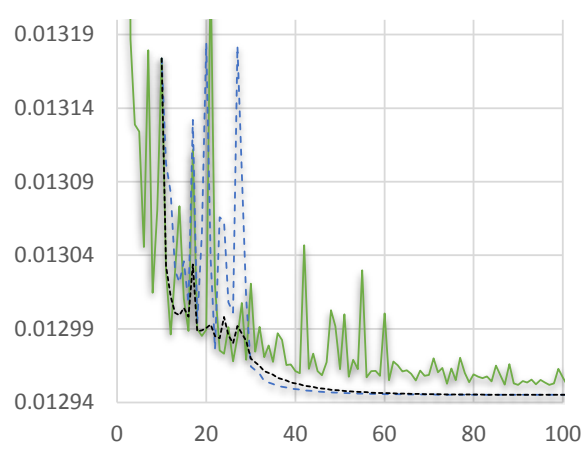

(k) mnist, $\lambda = 10^{-5}$

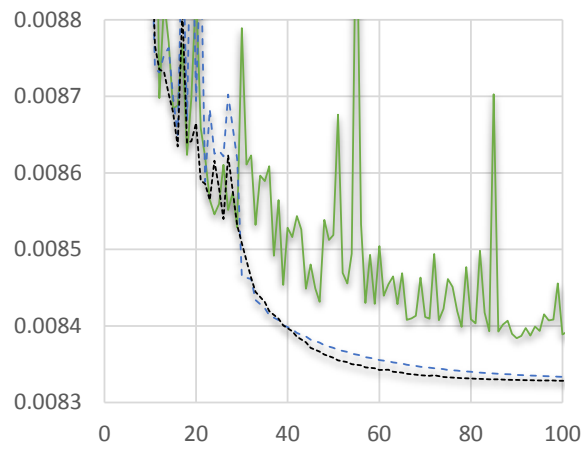

(l) mnist, $\lambda = 10^{-6}$

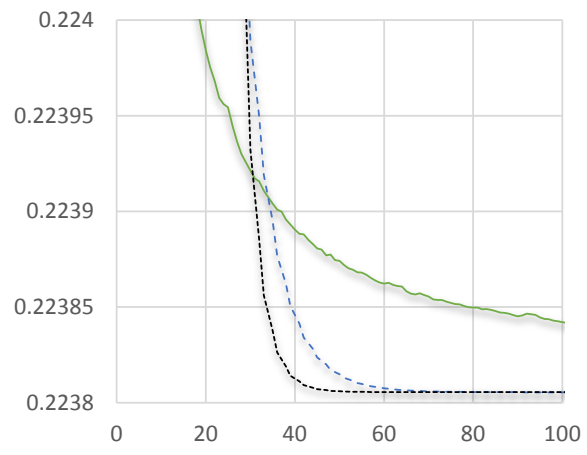

(m) rcv1, $\lambda = 10^{-3}$

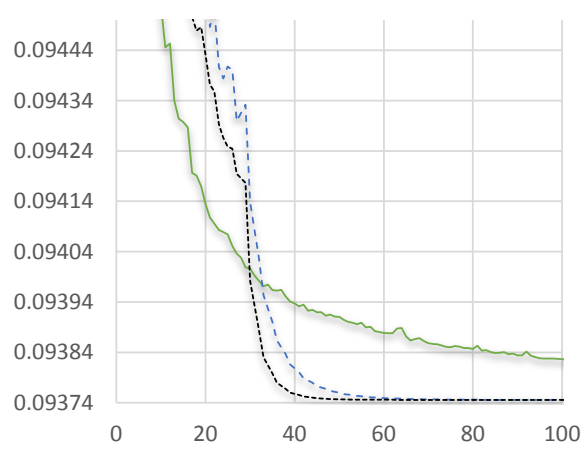

(n) rcv1, $\lambda = 10^{-4}$

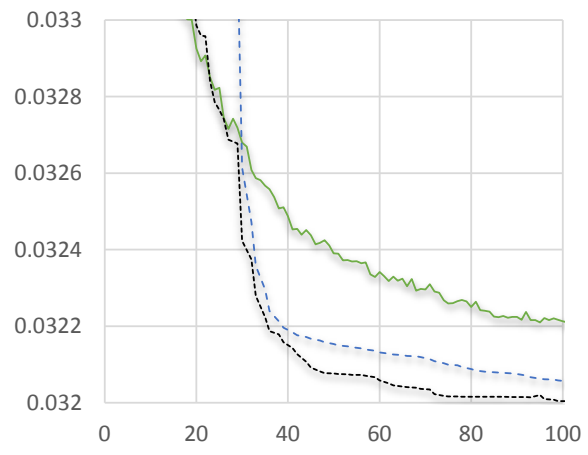

(o) rcv1, $\lambda = 10^{-5}$

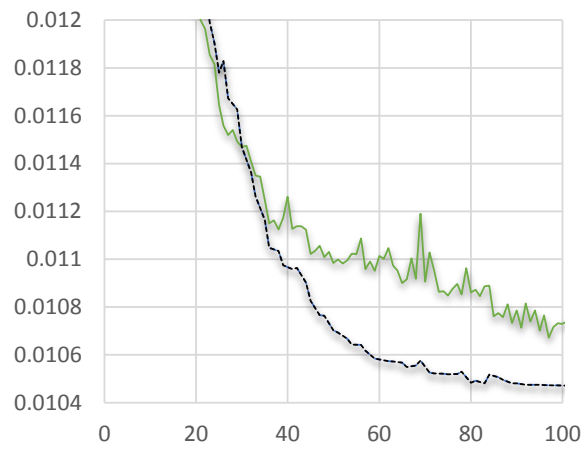

(p) rcv1, $\lambda = 10^{-6}$

Figure 5: Training error comparison between SGD and SVRG on $\ell_2$-regularized ERM with sigmoid loss. The best-tuned SGD is presented in solid green, the best-tuned SVRG with constant step length is presented in dashed blue, and the best-tuned SVRG is presented in **doted black**.

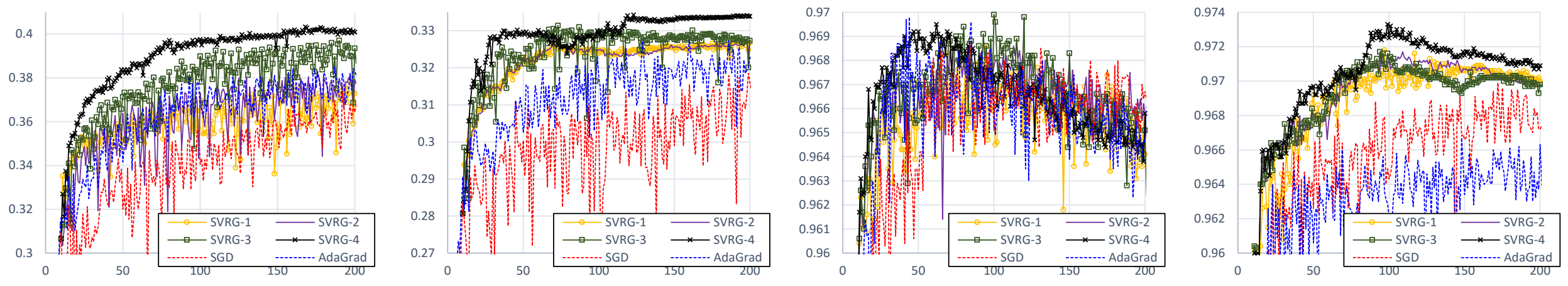

(a) cifar-10, $\lambda = 0$ (b) cifar-10, $\lambda = 10^{-3}$ (c) mnist, $\lambda = 0$ (d) mnist, $\lambda = 10^{-4}$

Figure 6: Test Accuracy Comparison on neural nets.